\documentclass[10pt]{article}

\usepackage{amsmath}
\usepackage{amssymb}
\usepackage{latexsym}
\usepackage{amsthm}

\newcommand{\init}{\big\vert_{t = 0}}
\newcommand{\smallinit}{\vert_{t = 0}}
 \newcommand{\half}{\frac{1}{2}}
\newcommand{\abs}[1]{\left\vert #1 \right\vert}
\newcommand{\bigabs}[1]{\bigl\vert #1 \bigr\vert}

\newcommand{\norm}[1]{\left\Vert #1 \right\Vert}
\newcommand{\bignorm}[1]{\bigl\Vert #1 \bigr\Vert}

\newcommand{\Sobnorm}[2]{\norm{#1}_{H^{#2}}}

\newcommand{\Sobdotnorm}[2]{\norm{#1}_{\dot{H}^{#2}}}

\newcommand{\Lpnorm}[2]{\norm{#1}_{L^{#2}}}

\newcommand{\Lxpnorm}[2]{\norm{#1}_{L_{x}^{#2}}}
\newcommand{\bigLxpnorm}[2]{\bignorm{#1}_{L_{x}^{#2}}}

\newcommand{\twonorm}[2]{\norm{#1}_{L^2#2}}
\newcommand{\bigtwonorm}[2]{\bignorm{#1}_{L^2#2}}
\newcommand{\inftynorm}[2]{\norm{#1}_{L^\infty#2}}

\newcommand{\mixednorm}[3]{\norm{#1}_{L_{t}^{#2}L_{x}^{#3}}}
\newcommand{\mixednormlocal}[3]{\norm{#1}_{L_{t}^{#2}L_{x}^{#3}(S_{T})}}

\newcommand{\LpHslocal}[3]{\norm{#1}_{L_{t}^{#2}H^{#3}(S_{T})}}
\newcommand{\LpdotHslocal}[3]{\norm{#1}_{L_{t}^{#2}\dot H^{#3}(S_{T})}}
\newcommand{\mixed}[2]{L_{t}^{#1}L_{x}^{#2}}

\newcommand{\bigmixednormlocal}[3]{\bignorm{#1}_{L_{t}^{#2}L_{x}^{#3}(S_{T})}}
\newcommand{\energy}[1]{\norm{#1}_{L_{t}^{\infty}H^{1}}}

\newcommand{\energylocal}[1]{\norm{#1}_{L_{t}^{\infty}H^{1}(S_{T})}}

\newcommand{\dotenergylocal}[1]{\norm{#1}_{L_{t}^{\infty}\dot H^{1}(S_{T})}}

\newcommand{\A}{\mathbf{A}}
\newcommand{\bfa}{\mathbf{a}}
\newcommand{\B}{\mathbf{B}}
\newcommand{\C}{\mathbb{C}}
\newcommand{\D}{\abs{D_{x}}}
\newcommand{\E}{\mathbf{E}}

\newcommand{\Proj}{\mathcal{P}}
\newcommand{\R}{\mathbb{R}}

\newcommand{\Fourier}{\mathcal{F}}

\newcommand{\innerprod}[2]{\left\langle \, #1 , #2 \, \right\rangle}

\DeclareMathOperator{\diag}{diag}
\DeclareMathOperator{\dv}{div}

\newtheorem{theorem}{Theorem}
\newtheorem{proposition}{Proposition}
\newtheorem{lemma}{Lemma}
\newtheorem*{corollary}{Corollary}

\theoremstyle{definition}

\newtheorem*{notation}{Notation}

\theoremstyle{remark}

\newtheorem*{plainremark}{Remark}
\newtheorem*{plainremarks}{Remarks}

\title{Nonrelativistic limit of Klein-Gordon-Maxwell to 
Schr\"odinger-Poisson}

\author{Philippe Bechouche, Norbert Mauser and Sigmund Selberg\\
Wolfgang Pauli Institut,\\
c/o Inst. f. Mathematik, \\
Universit\"at Wien, \\
Strudlhofgasse 4, A-1090 Wien}

\date{}

\begin{document}

\maketitle

\begin{abstract}
  We prove that in the nonrelativistic limit $c \to \infty$, where $c$ 
  is the speed of light, solutions of the Klein-Gordon-Maxwell 
  system on $\R^{1+3}$ converge in the energy space $C([0,T];H^{1})$ to
  solutions of a Schr\"odinger-Poisson system, under appropriate 
  conditions on the initial data. This requires the splitting of the 
  scalar Klein-Gordon field into a sum of two fields, corresponding, 
  in the physical interpretation, to electrons and positrons.
\end{abstract}

\section{Introduction}

\subsection{Klein-Gordon-Maxwell on $\R^{1+3}$}

The Klein-Gordon-Maxwell (abbreviated KGM) system on $\R^{1+3}$ reads
\begin{subequations}\label{realKGM}
\begin{align}
 \label{realKGMb}
  D_{\mu} D^{\mu} \phi &= c^{2} \phi,
  \\
  \label{realKGMa}
  \partial^{\nu} F_{\mu\nu}
  &= \tfrac{1}{c} \Im \left( \phi \overline{D_{\mu} \phi} \right).
\end{align}
\end{subequations}
In this paper we shall rely on the Coulomb gauge condition
\begin{equation}\label{CG}
  \dv \A = \partial^{i} A_{i} = 0,
\end{equation}
which has certain advantages for KGM, as demonstrated in \cite{KM}.

Here we use relativistic coordinates $x^{0} = ct \in \R, x = ( x^{1}, 
x^{2}, x^{3} ) \in \R^{3}$, where $c$ is the light speed.
Indices are raised and lowered relative to the
Minkowski metric with signature $-1,1,1,1$. The Einstein summation convention 
is in effect: Greek indices are summed over $0,1,2,3$, roman indices 
over $1,2,3$. We write $\partial_{\mu} = \tfrac{\partial}{\partial x^{\mu}}$.
Thus, $\partial_{0} = \tfrac{1}{c} \partial_{t}$, where
$\partial_{t} = \tfrac{\partial}{\partial t}$.
$\phi \in \C$ represents a particle field
and $F_{\mu\nu}$ is the electromagnetic field tensor, given 
in terms of a real potential $A_{\mu}$
by
\begin{equation}\label{EMtensor}
  F_{\mu \nu} = \partial_{\mu} A_{\nu} - \partial_{\nu} A_{\mu}.
\end{equation}
We split $A_{\mu}$ into its temporal part $A_{0}$ and its spatial 
part $\A = (A_{1},A_{2},A_{3})$. $D_{\mu}$ is the covariant derivative
$$
  D_{\mu} \phi = \partial_{\mu}\phi + \tfrac{i}{c} A_{\mu}\phi.
$$
Thus, since $\partial_{0} = \tfrac{1}{c} \partial_{t}$,
\begin{equation}\label{Ds}
  D_{0}\phi = \tfrac{1}{c}( \partial_{t} \phi + i A_{0} \phi),
  \qquad
  D_{j} \phi = \partial_{j} \phi + \tfrac{i}{c} A_{j} \phi \quad (j = 
  1,2,3).
\end{equation}
For $z \in \C$, $\Re z$ denotes the real part and $\Im z$ the imaginary part.
We also write $\nabla = (\partial_{1},\partial_{2},\partial_{3})$,
$\Delta = \partial_{i} \partial^{i} =
\partial_{1}^{2} + \partial_{2}^{2} + \partial_{3}^{2}$ and
$\square = \partial_{\mu} \partial^{\mu}
= - \tfrac{1}{c^{2}} \partial_{t}^{2} + \Delta$. 

For the convenience of the reader, let us briefly recall the derivation of \eqref{realKGM}.
First, write Maxwell's equations in the form
\begin{equation}\label{Maxwell}
  \partial^{\nu} F_{\mu\nu} = \tfrac{4\pi}{c} j_{\mu},
\end{equation}
where $j_{\mu}$ is the four-current density.
To translate \eqref{Maxwell} into classical notation, define the electric and 
magnetic field vectors by
\begin{equation}\label{EandB}
  \E = \nabla A_{0} - \tfrac{1}{c} \partial_{t} \A
  \quad \text{and} \quad \B = \nabla \times \A
\end{equation}
respectively. Then $F_{i0} = E_{i}$ and $F_{ij} = \epsilon_{ijk} B^{k}$, so 
\eqref{Maxwell} becomes
$$
  \dv \E = 4\pi \rho, \quad
  \nabla \times \B - \tfrac{1}{c} \partial_{t} \E = \tfrac{4\pi}{c} 
  \mathbf j,
$$
where $\rho = \tfrac{1}{c} j^{0}$ and $\mathbf j = 
(j^{1},j^{2},j^{3})$. The equations
$$
  \dv \B = 0, \quad \nabla \times \E + \tfrac{1}{c} \partial_{t} \B = 
  0,
$$
which follow from the definitions of $\E$ and $\B$, complete the Maxwell system in 
standard form.

Next, recall the free Klein-Gordon equation for a particle with rest mass $m$,
\begin{equation}\label{FreeKG}
  \square\phi = m^{2} c^{2} \phi.
\end{equation}
The associated current density
\begin{equation}\label{FreeKGcurrent}
  j_{\mu} = \tfrac{1}{2mi} \left( \phi \overline{\partial_{\mu} \phi}
  - \overline{\phi} \partial_{\mu} \phi \right)
  = \tfrac{1}{m} \Im \left( \phi \overline{\partial_{\mu} \phi} \right)
\end{equation}
satisfies the conservation law $\partial^{\mu} j_{\mu} = 0$. In terms of
$\rho = \tfrac{1}{c} j^{0}$ and $\mathbf j = (j^{1},j^{2},j^{3})$, 
this reads $\partial_{t} \rho + \dv \mathbf j = 0$.

The coupling of the free Klein-Gordon equation
to an electromagnetic field represented by $A_{\mu}$
is achieved by the so-called minimal substitution
$$
  \partial_{\mu} \longrightarrow D_{\mu}.
$$
Thus, setting $m = 1$ from now on, \eqref{FreeKG} transforms to 
\eqref{realKGMb}, and the current density \eqref{FreeKGcurrent}
transforms to
\begin{equation}\label{KGcurrent}
  j_{\mu} = \Im \left( \phi \overline{D_{\mu} \phi} \right)
  = \Im \left( \phi \overline{\partial_{\mu} \phi} \right)
  - \tfrac{1}{c} A_{\mu} \abs{\phi}^{2},
\end{equation}
which again satisfies\footnote{This follows from \eqref{realKGMb}, in view of
the identity
$\partial^{\mu} j_{\mu} = \Im \left( \phi \overline{D_{\mu} D^{\mu} \phi} \right)$.}
$\partial^{\mu} j_{\mu} = 0$. 
Substituting \eqref{KGcurrent} into Maxwell's equation \eqref{Maxwell} and 
dropping the factor $4\pi$ gives \eqref{realKGMa}.

The system \eqref{realKGM} can also be derived from Hamilton's 
principle using the Lagrangian density
$$
  \mathcal L = - \frac{1}{4} F_{\mu\nu} F^{\mu\nu} - \half \left( 
  D_{\mu} \phi \overline{ D^{\mu} \phi } + c^{2} \phi \overline{\phi} 
  \right).
$$
Since $\mathcal L$ does not depend explicitly on $x^{\mu}$, the 
energy-momentum tensor
$$
  T_{\mu\nu} = \frac{\partial \mathcal L}{\partial (\partial^{\nu} 
  A_{\lambda})} \partial_{\mu} A_{\lambda}
  + \frac{\partial \mathcal L}{\partial (\partial^{\nu} 
  \phi)} \partial_{\mu} \phi
  + \frac{\partial \mathcal L}{\partial (\overline{\partial^{\nu} 
  \phi}) } \overline{\partial_{\mu} \phi}
  - \mathcal L \delta_{\mu\nu}
$$
satisfies $\partial^{\nu} T_{\mu\nu} = 0$. See, for example, 
\cite[Chapter 12]{Goldstein}. This tensor turns out not to be 
symmetric, but we can symmetrize it by the same trick that one uses 
for the Maxwell Lagrangian (see \cite[pp 583--584]{Goldstein}). Thus, 
we set $T'_{\mu\nu} = T_{\mu\nu} - \partial^{\lambda} \left(
F_{\lambda\nu} A_{\mu} \right)$. Then we still have the conservation 
law $\partial^{\nu} T'_{\mu\nu} = 0$, which in particular implies
\begin{equation}\label{ConsLaw}
  \mathcal E(t) = \int_{\R^{3}} T'_{00}(t,x) \, dx = \text{const.}
\end{equation}
A calculation reveals that
\begin{equation}\label{T00}
  T'_{00} = \half \left( \abs{D_{0}\phi}^{2} + \sum\nolimits_{1}^{3} 
  \abs{D_{i} \phi}^{2} + c^{2} \abs{\phi}^{2} 
  + \E^{2} + \B^{2} \right),
\end{equation}
and it is a fortunate fact that this density is non-negative.

\subsection{Main result}\label{MainResult}

The system \eqref{realKGM}, \eqref{CG} can be reformulated as follows:\footnote{Cf. 
\cite[Eqs. (1.7a--c)]{KM}. There, however, the light
speed $c = 1$ and the rest mass is zero, so that $M^{2} = - \Delta$
in \eqref{KGMc}.}
\begin{subequations}\label{KGM}
\begin{align}
  \label{KGMc}
  \left( \partial_{t}^{2} + M^{2} \right) \phi &= 2ic \A \cdot 
  \nabla \phi - 2i A_{0} \partial_{t} \phi - i (\partial_{t} A_{0}) 
  \phi + (A_{0}^{2} - \A^{2}) \phi,
  \\
  \label{KGMa}
  \Delta A_{0} &= - \tfrac{1}{c^{2}} \Im \left( \phi \overline{ 
  \partial_{t} \phi} \right) + \tfrac{1}{c^{2}} \abs{\phi}^{2} A_{0},
  \\
  \label{KGMb}
  \square \A &= - \tfrac{1}{c} \Proj \left( \Im \left( \phi \overline{ 
  \nabla \phi} \right) \right) + \tfrac{1}{c^{2}} \Proj \bigl( 
  \abs{\phi}^{2} \A \bigr),
\end{align}
\end{subequations}
where
\begin{equation}\label{Moperator}
  M = M(c) = \sqrt{c^{4} - c^{2} \Delta}
\end{equation}
and $\Proj$ is the projection onto the divergence free vector fields on $\R^{3}$.
In terms of the Riesz operators $R_{i} = (-\Delta)^{-1/2} \partial_{i}$, 
\begin{equation}\label{ProjDef}
  \Proj X^{i} = X^{i} - \sum_{j = 1}^{3} R_{i} R_{j} X^{j} \quad (i = 
  1,2,3).
\end{equation}
In fact, expanding \eqref{realKGMa} using \eqref{EMtensor} and 
\eqref{CG}, one obtains \eqref{KGMa} for $\mu = 0$, and for $\mu = 
1,2,3$ one gets
$$
  \square \A- \partial_{0} \nabla A_{0} =
  - \tfrac{1}{c} \Im \left( \phi (\overline{D_{i} \phi})_{i = 1,2,3} \right).
$$
Apply $\Proj$ to both sides, use the identity $\Proj \nabla \equiv 
0$, and observe that \eqref{CG} implies $\Proj \A = \A$, to
obtain \eqref{KGMb}. Finally, \eqref{realKGMb} expands to give \eqref{KGMc},
if we use \eqref{CG}.

We specify finite energy initial data at time $t = 0$:\footnote{Equivalently,
instead of \eqref{KGMdata1} we could specify data for $\E$ and $\B$ in $L^{2}$.}
\begin{subequations}\label{KGMdata}
\begin{alignat}{2}
  \label{KGMdata1}
  \A \init &= \bfa_{0}(c) \in \dot H^{1}, &\quad \partial_{t} \A \init &= 
  \bfa_{1}(c) \in L^{2},
  \\
  \label{KGMdata2}
  \phi \init &= \phi_{0}(c) \in H^{1}, &\quad \partial_{t} \phi \init &= 
  \phi_{1}(c) \in L^{2}.
\end{alignat}
\end{subequations}
Here $H^{s} = H^{s}(\R^{3})$ is the Sobolev space with norm 
$\Sobnorm{f}{s} = \bigtwonorm{(1 + \abs{\xi}^{2})^{s/2} \widehat 
f}{_{\xi}}$, where $\widehat f(\xi)$ is the Fourier transform of $f(x)$, and $\dot H^{s}$
denotes the corresponding homogeneous space, with
norm $\Sobdotnorm{f}{s} = \bigtwonorm{\abs{\xi}^{s} \widehat f}{_{\xi}}$.

In view of the Coulomb condition \eqref{CG}, we must assume
\begin{equation}\label{DivFreeData}
  \dv \bfa_{0} = \dv \bfa_{1} = 0.
\end{equation}
Then \eqref{CG} is implicit in the system \eqref{KGM},
since \eqref{KGMb} implies $\square (\dv \A) = 0$, so that 
\eqref{DivFreeData} persists in time.

Klainerman and Machedon \cite{KM} proved that \eqref{KGM} is locally
well-posed\footnote{In \cite{KM} the rest mass is assumed to be
zero, but it is a trivial matter to modify the 
proof of local well-posedness to handle the linear term introduced when the
mass is positive. For the convenience of the reader, we give the details
in an appendix.} for initial data \eqref{KGMdata},
and then by conservation of energy they 
obtained global well-posedness for such data.
Using linear Strichartz estimates for the homogeneous wave equation 
$\square u = 0$, it is possible to prove local well-posedness if one 
assumes slightly more regularity of the data.
To get the result proved in \cite{KM}, however, requires certain bilinear 
generalizations of Strichartz' $L^{4}$ estimate (see \cite[Section 2]{KM})
to handle the first terms on the right hand sides of (\ref{KGM}b,c). 
A key point is that, due to the Coulomb gauge condition, these terms have
a null form structure, without which the estimates would in fact fail.
Here we will need modifications of these estimates (see section 
\ref{SpacetimeEstimates}) where the wave 
operator $\square$ may be replaced by $i\partial_{t} \pm (M-c^{2})$.
The latter essentially behaves like the Schr\"odinger operator at
frequency $\lesssim c$, and like the wave operator at frequency $\gg c$.

The global solutions of \eqref{KGM}, \eqref{KGMdata} obtained in 
\cite{KM} have the regularity
\begin{equation}\label{KMsolutionReg1}
  \partial_{\mu} A_{\nu} \in C(\R;L^{2}),
  \quad \phi \in C(\R;H^{1}) \cap C^{1}(\R;L^{2}).
\end{equation}
Moreover, for every $0 < T < \infty$ (see the Main Theorem and Propositions 3.2
and 2.3 in \cite{KM})
\begin{gather}\label{KMsolutionReg2}
  \int_{0}^{T} \left( \twonorm{\square \A(t)}{} + \twonorm{\square 
  \phi(t)}{} \right) \, dt < \infty,
  \\
  \label{KMsolutionReg3}
  \int_{0}^{T} \left( \Lpnorm{\nabla A_{0}(t)}{3}
  + \inftynorm{A_{0}(t)}{} \right) \, dt < \infty.
\end{gather}

The question considered in this paper is what happens to the solutions 
as $c \to \infty$. Let us first state our main result, and then in the next section
we motivate it.

Throughout the paper, the $O,o$ notation refers to the limit $c \to \infty$.
The following notation is used for function spaces.
If $X$ is a Banach space of functions on $\R_{x}^{3}$, we denote by 
$L_{t}^{p} X$ the space with norm $\norm{u}_{L_{t}^{p}X} = ( 
\int_{-\infty}^{\infty} \norm{u(t,\cdot)}_{X}^{p} \, dt )^{1/p}$, with the 
usual modification if $p = \infty$. The localization of this norm 
to $S_{T} = [0,T] \times \R^{3}$ is denoted 
$\norm{u}_{L_{t}^{p}X(S_{T})}$.

\begin{theorem}\label{Thm2}
Suppose $(A_{0},\A,\phi)$ solve \eqref{KGM} on $\R^{1+3}$ with data 
\eqref{KGMdata} such that \footnote{This is equivalent to 
$\twonorm{\E(t = 0)}{} + \twonorm{\B(t = 0)}{} = O(1)$.}
$$
  \Sobdotnorm{\bfa_{0}(c)}{1} + \tfrac{1}{c} 
  \twonorm{\bfa_{1}(c)}{} = O(1),
$$
and such that the limits
\begin{equation}\label{DataLimits}
  \alpha = \lim_{c \to \infty} \phi_{0}(c) \quad \text{and} \quad
  \beta = \lim_{c \to \infty} M^{-1} \phi_{1}(c)
\end{equation}
exist in $H^{1}$. Split $\phi = \phi^{+} + \phi^{-}$, where
\begin{equation}\label{FreeSplitting}
  \phi^{\pm} = \half \left\{ \phi \pm i M^{-1} \partial_{t} \phi \right\},
\end{equation}
and set
\begin{equation}\label{PsiDef}
  \psi^{\pm} = e^{\pm itc^{2}} \phi^{\pm}.
\end{equation}
Let $(u,v^{+},v^{-})$ be the solution 
of the Schr\"odinger-Poisson system
\begin{subequations}\label{PS}
  \begin{gather}
    \label{PSb}
    \Bigl( i \partial_{t} \pm \frac{\Delta}{2} \Bigr) v^{\pm} = u v^{\pm},
    \\
     \label{PSa}
    \Delta u = - \bigabs{v^{+}}^{2} + \bigabs{v^{-}}^{2},
  \end{gather}
\end{subequations}
on $\R^{1+3}$, with initial data $v^{\pm} \init = \tfrac{1}{2}(\alpha \pm i \beta)$.
Then for every $0 < T < \infty$,
\begin{align}
  \label{Convergence}
  \energylocal{\psi^{\pm} - v^{\pm}} &\longrightarrow 0,
  \\
  \label{A02uConvergence}
  \dotenergylocal{A_{0} - u} &\longrightarrow 0,
\end{align}
as $c \to \infty$.
\end{theorem}

\begin{plainremarks} (i) From \eqref{Convergence} it follows that
$\phi - e^{-itc^{2}} v^{+} - e^{+itc^{2}} v^{-} \to 0$ in 
$L_{t}^{\infty}H^{1}(S_{T})$, since $\phi = \phi^{+} + \phi^{-} =
e^{-itc^{2}} \psi^{+} - e^{+itc^{2}} \psi^{-}$.

\medskip\noindent
(ii) We will in fact prove that, for $1 \le r \le 3/2$,
\begin{equation}\label{DeltaA0minusu}
  \mixednormlocal{\Delta (A_{0}-u)}{\infty}{r} \longrightarrow 0
  \quad \text{as} \quad c \longrightarrow \infty.
\end{equation}
Then by Sobolev embedding and the fact that (see \cite[Proposition III.3]{St})
\begin{equation}\label{RieszCorollary}
  \Lpnorm{\partial_{i} \partial_{j} f}{p} \lesssim \Lpnorm{\Delta f}{p}
  \quad \text{for} \quad  1 < p < \infty,
\end{equation}
we get
$\mixednormlocal{A_{0}-u}{\infty}{r} \to 0$ for $3 < r < \infty$
and $\mixednormlocal{\nabla (A_{0}-u)}{\infty}{r} \to 0$ for $3/2 < r \le 3$.
In particular, this gives \eqref{A02uConvergence}.

\medskip\noindent
(iii) The system \eqref{PS} is globally well-posed in $L^{2}$, as proved by
Castella \cite{Castella}. See section \ref{PSbounds} for details.
\end{plainremarks}

\medskip

In the literature there are some results concerning the nonrelativistic limit of
the linear Klein-Gordon equation with a fixed electromagnetic potential, see
\cite{CC} and \cite{Ve}, but there are no previous results for the coupled
Klein-Gordon-Maxwell system. Moreover, these papers only treat the static case,
i.e., the potential is time-independent. The nonrelativistic limit for the related Dirac
equation with time-dependent external potential was treated in \cite{BMP}.

However, we have been made aware of recent, completely independent work
of Masmoudi and Nakanishi \cite{MaNa}, who have obtained results similar
to ours.

To motivate our result, we consider a simple but instructive example, 
namely the free Klein-Gordon equation.

\subsection{Model case: Free Klein-Gordon}\label{ModelCase}

In the absence of an electromagnetic field, \eqref{KGMc} would reduce to 
the free Klein-Gordon equation
\begin{equation}\label{FreeKGsection1}
  \square \phi = c^{2} \phi.
\end{equation}
The latter can be derived as a relativistic analogue of the free Schr\"odinger equation.
Indeed, recall the quantum mechanical principle whereby
classical quantities are replaced by operators:\footnote{We set 
Planck's constant equal to $1$.}
\begin{itemize}
  \item[] Energy $\qquad E \longrightarrow i \frac{\partial}{\partial 
  t}$,
  \item[] Momentum $\qquad \mathbf p \longrightarrow \frac{1}{i} \nabla$.
\end{itemize}
Thus, from the nonrelativistic energy of a free particle with unit rest 
mass,
$$
  E = \frac{\mathbf p^{2}}{2},
$$
one obtains the free Schr\"odinger equation $i \partial_{t} \psi = - 
\frac{\Delta}{2} \psi$. 
Proceeding instead from the relativistic energy-momentum relation
\begin{equation}\label{relEM}
  E = \sqrt{ c^{2} \mathbf p^{2} + c^{4} }
\end{equation}  
gives, in our notation,
$i\partial_{t} \phi = M \phi$.
Squaring this gives $- \partial_{t}^{2} \phi = M^{2} \phi$, which we
can write in the compact and obviously Lorentz invariant form 
\eqref{FreeKGsection1}.

We now ask, conversely, whether Klein-Gordon reverts to
Schr\"odinger in the non-relativistic limit $c \to \infty$.
Reversing the above steps, we see that we first have to formally take 
square roots of \eqref{FreeKGsection1}, written in the form $- \partial_{t}^{2} 
\phi = M^{2} \phi$. Then there will naturally be two separate fields
$\phi^{\pm}$, corresponding to positive and negative square roots, and
solving
\begin{equation}\label{SqrtOfKG}
  i \partial_{t} \phi^{\pm} = \pm M \phi^{\pm}.
\end{equation}
For $\phi^+$, this corresponds to the relation \eqref{relEM}, but before we 
can let $c \to \infty$ there, we clearly have to subtract the dominant term 
$c^{2}$, the rest energy. Thus, we note that if $E$ is given by 
\eqref{relEM}, then
\begin{equation}\label{Eminusc2}
  E - c^{2} = \frac{\mathbf p^{2}}{1 + E^{2}/c^{2}} \to \frac{\mathbf p^{2}}{2}
\end{equation}
as $c \to \infty$. In the case of negative energy, i.e., for $\phi^-$,
we have to \emph{add} the rest energy, of course.  Since
the kinetic energy $E$ corresponds to the Fourier
variable of $-t$, while $\mathbf p$ corresponds to that of $x$, this 
procedure of subtracting or adding the rest energy corresponds in 
physical space to multiplication by the oscillating factors $e^{+ itc^{2}}$
or $e^{- itc^{2}}$, respectively. Since $\phi^{+}$ (resp. $\phi^{-}$)
corresponds to positive (resp. negative) energy states, we can think
of it as representing electrons (resp. positrons).

The above heuristics suggest that to analyze the limit $c \to \infty$, 
a solution $\phi$ of \eqref{FreeKGsection1} must be split
$\phi = \phi^{+} + \phi^{-}$,
where $\phi^{\pm}$ solve \eqref{SqrtOfKG}, hence \eqref{FreeKGsection1}.
But this forces the initial constraints
$$
  \phi_{0} = \phi_{0}^{+} + \phi_{0}^{-},
  \quad
  i\phi_{1} = M \phi_{0}^{+} - M\phi_{0}^{-},
$$
where $\phi_{0} = \phi \init$, $\phi_{1} = \partial_{t} \phi \init$ 
and $\phi_{0}^{\pm} = \phi^{\pm} \init$.
Solving for $\phi_{0}^{\pm}$ gives
\begin{equation}\label{FreeInitialSplitting}
   \phi_{0}^{\pm} = \half \left\{ \phi_{0} \pm i M^{-1} 
   \phi_{1} \right\}.
 \end{equation}
Next, we either subtract (for the electron) or add (for the positron)
the rest energy. That is, we define $\psi^{\pm}$ by \eqref{PsiDef}.
Then by \eqref{SqrtOfKG} and \eqref{FreeInitialSplitting},
$$
  i \partial_{t} \psi^{\pm} = \pm ( M - c^{2} ) \psi^{\pm},
  \quad \psi^{\pm} \init = \half \left\{ \phi_{0} \pm i M^{-1} 
   \phi_{1} \right\}.
$$
Observe that the limits $\lim_{c \to \infty} \psi^{\pm} \init$ exist in some 
space if and only if the limits in \eqref{DataLimits}
exist, and since the Fourier symbol of $M - c^{2}$ is (cf. \eqref{Eminusc2})
\begin{equation}\label{Mminusc2Symbol}
  \sqrt{ c^{2} \abs{\xi}^{2} + c^{4}} - c^{2}
  = \frac{\abs{\xi}^{2}}{1 + \sqrt{1 + \abs{\xi}^{2}/c^{2}}}
  \longrightarrow \frac{\abs{\xi}^{2}}{2} \quad \text{as} \quad c \to 
  \infty,
\end{equation}
we would then expect $\psi^{\pm}$ to converge to the solutions of
\begin{equation}\label{LimitEq}
  i \partial_{t} v^{\pm} = \mp \frac{\Delta}{2} v^{\pm},
  \quad v^{\pm} \init = \half( \alpha \pm i \beta ).
\end{equation}
This is of course easy to verify directly here, since 
$\psi^{\pm}= e^{\mp i t (M - c^{2})}\phi_{0}^{\pm}$.
So if the limits \eqref{DataLimits} exist in $H^{s}$, say,
then it follows by the dominated convergence theorem that $\psi^{\pm}$ converges
in $C([0,T]; H^{s})$ to $v^{\pm}(t) = e^{\pm i t \Delta/2} (\alpha \pm 
i \beta)/2$, which solves \eqref{LimitEq}.

\begin{plainremarks}
(i) In this example we defined $\phi^{\pm}$ as the 
solutions of \eqref{SqrtOfKG} with data \eqref{FreeInitialSplitting}, but 
this is equivalent to using \eqref{FreeSplitting} at each time $t$. 
Indeed, taking a time derivative in \eqref{FreeSplitting} and using the fact 
that $\phi$ solves \eqref{FreeKGsection1}, one obtains 
\eqref{SqrtOfKG}.

\medskip
\noindent
{(ii)} The splitting \eqref{FreeSplitting} of the Klein-Gordon field
into an ``electron'' part and a ``positron'' part corresponds exactly
to the splitting of the 4-spinor of the Dirac equation used in \cite{BMP}.
To see this, write \eqref{FreeKGsection1} as a first order system
$$
 \partial_t \begin{pmatrix} \phi^{(0)} \\ \phi^{(1)} \end{pmatrix}
  = \begin{pmatrix} 0 & I \\ - M^2 & 0 \end{pmatrix}
  \begin{pmatrix} \phi^{(0)} \\ \phi^{(1)} \end{pmatrix},
$$
where $\phi^{(0)} = \phi$ and $\phi^{(1)} = \partial_t \phi$.
The $2 \times 2$ matrix on the right hand side has two eigenvalues,
$\lambda^\pm = \mp iM$. The eigenspace projections corresponding
to $\lambda^\pm$ are, respectively,
$$
  \Pi_\pm = 
  \half \begin{pmatrix} I & \pm iM^{-1} \\ \mp iM & I \end{pmatrix},
$$
exactly as in \cite[Eq. (1.21)]{BMP}.
\end{plainremarks}

\subsection{A priori bounds for KGM}\label{MainResults}

The previous example shows that the splitting $\phi = \phi^{+} + \phi^{-}$
defined by \eqref{FreeSplitting}, and used in Theorem \ref{Thm2},
is motivated by the free Klein-Gordon equation. Recall that the coupling of the
latter to the electromagnetic field $A_{\mu}$ is achieved by the minimal
substitution $\partial_{\mu} \to D_{\mu}$, which 
transforms \eqref{FreeSplitting} to
\begin{equation}\label{PhiSplitting}
  \phi^{\pm} = \half \left\{ \phi \pm M^{-1} ( i \partial_{t} \phi - 
  A_{0} \phi) \right\}.
\end{equation}
Since it turns out (see Theorem \ref{Thm1}) that
$$
  \energylocal{M^{-1}(A_{0}\phi)} = O(c^{-1})
$$
under the hypotheses of Theorem \ref{Thm2}, it is clear that as far 
as Theorem \ref{Thm2} is concerned, it is immaterial whether we use 
\eqref{FreeSplitting} or \eqref{PhiSplitting}.
The latter, however, is more natural to work with, since the evolution
equations satisfied by $\psi^{\pm}$ turn out to be much nicer. (In particular, 
if one considers not only convergence of $\psi^{\pm} \to v^{\pm}$ but
also $\partial_{t} \psi^{\pm} \to \partial_{t} v^{\pm}$, then 
\eqref{PhiSplitting} must be used.)
In fact, in section \ref{AprioriBounds} we prove the following.

\begin{lemma}\label{SplitSystemLemma}
If $(A_{0},\A,\phi)$ solves \eqref{KGM}, $\phi^{\pm}$ is defined by 
\eqref{PhiSplitting} and $\psi^{\pm} = e^{\pm itc^{2}} \phi^{\pm}$, 
then
\begin{subequations}\label{KGM'}
\begin{align}
  \label{KGM'c}
  L^{\pm} \psi^{\pm} &= A_{0} \psi^{\pm} \pm e^{\pm itc^{2}} R,
  \\
  \label{KGM'a}
  \Delta A_{0} &= - \Re \left( (\phi^{+} + \phi^{-}) \tfrac{M}{c^{2}} 
  \bigl(\overline{\phi^{+}} - \overline{\phi^{-}} \bigr) \right),
  \\
  \label{KGM'b}
  \square \A &= - \tfrac{1}{c} \Proj \left( \Im \left( \phi \overline{ 
  \nabla \phi} \right) \right) + \tfrac{1}{c^{2}} \Proj \bigl( 
  \abs{\phi}^{2} \A \bigr),
\end{align}
\end{subequations}
where
\begin{gather}
  \label{Loperator}
  L^{\pm} = L^{\pm}(c) = i \partial_{t} \mp (M-c^{2}),
  \\
  \label{Remainder}
  R = \half M^{-1} \left\{ - 2ic \A \cdot \nabla \phi + 
  [A_{0},M-c^{2}] (\phi^{+} - \phi^{-}) + \A^{2} \phi \right\}
\end{gather}
and $[A_{0},M-c^{2}]$ is the commutator:
\begin{equation}\label{Commutator}
  [A_{0},M-c^{2}] \phi = A_{0} (M-c^{2}) \phi - (M-c^{2}) (A_{0} \phi).
\end{equation}
\end{lemma}

Since $M-c^{2}$ behaves like $- \tfrac{\Delta}{2}$ as $c \to \infty$ 
(cf. \eqref{Mminusc2Symbol}), and since it turns out that $R$ 
vanishes in the limit (see Theorem \ref{Thm1}), it is not surprising 
that \eqref{KGM'c} tends to \eqref{PSb}. Similarly, to motivate the 
convergence of \eqref{KGM'a} to \eqref{PSa}, observe that
expansion of RHS\eqref{KGM'a} gives
\begin{equation}\label{DeltaA0Expansion}
  \Delta A_{0} = - \bigabs{\psi^{+}}^{2} + \bigabs{\psi^{-}}^{2} + 
  \tfrac{1}{c^{2}} R',
\end{equation}
where
\begin{equation}\label{A0Error}
\begin{split}
  R' &=
  - \Re \left( \psi^{+} (M-c^{2}) \overline{\psi^{+}} 
  \right)
  + \Re \left( \psi^{-} (M-c^{2}) \overline{\psi^{-}} 
  \right)
  \\
  &\quad
  + \Re \left( e^{-2itc^{2}} \psi^{+} (M-c^{2}) \overline{\psi^{-}} 
  \right)
  - \Re \left( e^{+2itc^{2}} \psi^{-} (M-c^{2}) \overline{\psi^{+}} 
  \right).
\end{split}
\end{equation}

The main difficulty in proving Theorem \ref{Thm2} is to obtain \emph{a 
priori} bounds as $c \to \infty$. The bounds obtained from the 
conservation of the KGM energy are not good enough.
For example, energy conservation gives $\energy{\phi} = O(c)$
(see section \ref{EnergyEstimates}), but this can be 
improved to $O(1)$ (on finite time intervals) using spacetime 
estimates of Strichartz type.
Energy conservation does, however, give the important global-in-time 
bound $\mixednorm{\phi^{\pm}}{\infty}{2} = O(1)$, which is not 
surprising in view of the fact that for the limiting system \eqref{PS},
the $L^{2}$ norms of $v^{\pm}$ are exactly conserved in time.

The main estimates are contained in the following theorem.

\begin{theorem}\label{Thm1}
  Suppose $(A_{0},\A,\phi)$ solve \eqref{KGM} with data 
  \eqref{KGMdata} satisfying
  \begin{gather}
    \label{AdataBound}
    \Sobdotnorm{\bfa_{0}(c)}{1} + \tfrac{1}{c} 
    \twonorm{\bfa_{1}(c)}{} = O(1),
    \\
    \label{PhidataBound}
    \Sobnorm{\phi_{0}(c)}{1} + \Sobnorm{M^{-1} \phi_{1}(c)}{1} = O(1).
  \end{gather}
  Then we have the global-in-time bound
  \begin{equation}\label{L2bound}
    \mixednorm{\phi^{\pm}}{\infty}{2} = 
    \mixednorm{\psi^{\pm}}{\infty}{2} = O(1).
  \end{equation}
  Moreover, for every $0 < T < \infty$,
  \begin{enumerate}
    \item\label{Thm1A}
    $\energylocal{\phi^{\pm}} = \energylocal{\psi^{\pm}} = O(1)$,
    \item\label{Thm1B}
    $\dotenergylocal{\A} + \tfrac{1}{c} \mixednormlocal{\partial_{t} 
    \A}{\infty}{2} = O(1)$,
    \item\label{Thm1C}
    $\LpHslocal{L^{\pm}\psi^{\pm}}{1}{1} + c \mixednormlocal{\square 
    \A}{1}{2} = O(1)$,
    \item\label{Thm1F}
    $\mixednormlocal{\nabla A_{0}}{\infty}{r} = O(1)$ for $3/2 < r \le 3$,
    \item\label{Thm1G}
    $\LpHslocal{R}{1}{1} = O(c^{-1/2})$,
    \item\label{Thm1H}
    $\energylocal{M^{-1}(A_{0}\phi)} = O(c^{-1})$,
  \end{enumerate}
  where $\phi^{\pm}$, $\psi^{\pm}$ and $R$ are given by 
  \eqref{PhiSplitting}, \eqref{PsiDef} and \eqref{Remainder}.
\end{theorem}

\begin{plainremark} Of course, \eqref{Thm1A} holds just as well for $\phi^{\pm}$ 
defined by \eqref{FreeSplitting}, in view of \eqref{Thm1H}.
The same remark applies to \eqref{L2bound}. In fact, from the 
proof of the latter, given in section \ref{EnergyEstimates}, we have 
the global-in-time bound $\mixednorm{M^{-1}(A_{0}\phi)}{\infty}{2} = O(c^{-1/2})$.
\end{plainremark}

The rest of this paper is organized as follows: The next section 
deals with the limit system \eqref{PS}, and in section 
\ref{Inequalities} we collect some inequalities that are used 
repeatedly. In section \ref{EnergyEstimates} we use energy 
conservation to prove \eqref{L2bound}, and section \ref{SpacetimeEstimates}
deals with linear and bilinear spacetime estimates for the operators $\square$ 
and $L^{\pm}$. In section \ref{AprioriBounds} we prove parts 
\eqref{Thm1A}--\eqref{Thm1H} of Theorem \ref{Thm1}, and finally in 
section \ref{H1Convergence} we prove the main result, Theorem 
\ref{Thm2}.

\begin{notation} Throughout the paper, the following conventions are 
in effect:
\begin{itemize}
  \item $\lesssim$ means $\le$ up to multiplication by an absolute, 
  positive constant. $X \sim Y$ stands for $X \lesssim Y \lesssim X$.
  \item The $O,o$ notation always refers to the limit $c \to \infty$.
  \item $K, \delta$ and $N$ denote absolute, positive constants which 
  may change from line to line. $\sigma(T)$ denotes the function 
  $K(T^{\delta} + T^{N})$ and $P(x)$ is the polynomial $x + x^{N}$.
  \item For exponents we use the standard shorthand $p^{+}$ (resp. 
  $p^{-}$) for $p + \varepsilon$ (resp. $p - \varepsilon$), where 
  $\varepsilon > 0$ is sufficiently small.
  See, e.g., Lemma \ref{EllipticLemma} in section \ref{Inequalities}.
  \item $\chi$ is a smooth cut-off on $\R^{3}$ such that $\chi(\xi) = 1$ for 
  $\abs{\xi} \le 1$ and $\chi(\xi) = 0$ for $\abs{\xi} \ge 2$. 
  Moreover, we assume that $\chi$ is radial, and we write $\chi(\xi)$ 
  and $\chi(r = \abs{\xi})$ interchangeably. We use $\chi(\xi/c)$ to 
  split functions $f(x)$ into low ($\lesssim c$) and high ($\gg c$) frequencies:
  \begin{equation}\label{LowHigh}
    f = f * \theta_{c} + f * (1-\theta_{c}) = f_{l} + f_{h},
  \end{equation}
  where $\theta_{c}$ is the inverse Fourier transform of $\chi(\xi/c)$.
  Then $\Lpnorm{\theta_{c}}{1}$ does not depend on $c$, so
  $\Lpnorm{f_{l}}{p}, \Lpnorm{f_{h}}{p} \lesssim \Lpnorm{f}{p}$ for $1 \le p \le 
  \infty$ by Young's inequality.
\end{itemize}
\end{notation}

\subsection{$H^{1}$ bounds for Schr\"odinger-Poisson}\label{PSbounds}

Global well-posedness in $L^2$ for the Schr\"odinger-Poisson system \eqref{PS}
follows from the work of Castella \cite{Castella}. In fact, since the
$L^{2}$ norms of $v^{\pm}$ are conserved:
\begin{equation}\label{L2conservation}
  \twonorm{v^{\pm}(t)}{} = \twonorm{v^{\pm}(0)}{} \quad \text{for}
  \quad t \ge 0,
\end{equation}
it is enough to prove local well-posedness for $L^{2}$ data.
It is then easy to obtain $L^{2}$ bounds for $\nabla v^{\pm}$ on
finite time intervals. For the convenience of the reader, and since
a similar but more involved argument will be used in the proof of Theorem \ref{Thm1}
(see section \ref{Zestimate}), we include here a short proof of these facts.
Thus, we prove:

\begin{lemma}\label{PSlemma}(Cf. \cite{Castella}.)
  The system \eqref{PS} is globally well-posed in $L^{2}$, and for 
  $H^{1}$ initial data we have
  \begin{equation}\label{PSH1bound}
    \energylocal{v^{\pm}} < \infty
  \end{equation}
  for all $T < \infty$.
\end{lemma}

So assume $(u,v^{+},v^{-})$ is a solution of \eqref{PS}, and let us 
derive some \emph{a priori} estimates for $v^{+}$ (the argument for 
$v^{-}$ is of course the same). Writing $\innerprod{f}{g} = 
\int_{\R^{3}} f \overline g \, dx$, we have
\begin{align*}
  &\frac{d}{dt} \half \innerprod{\nabla v^{+}}{\nabla v^{+}}
  = \Re \innerprod{\nabla \partial_{t} v^{+}}{\nabla v^{+}}
  \\
  &\quad = \Im \innerprod{-\tfrac{\Delta}{2} (\nabla v^{+})
  + u \nabla v^{+} + (\nabla u) v^{+}}{\nabla v^{+}}
  = \Im \innerprod{(\nabla u) v^{+}}{\nabla v^{+}},
\end{align*}
since $\Delta$ and $u$ are self-adjoint. But
\begin{align*}
  \Im \innerprod{(\nabla u) v^{+}}{\nabla v^{+}}
  &\le \Lxpnorm{\nabla u}{3} \Lxpnorm{v^{+}}{6} \Lxpnorm{\nabla v^{+}}{2}
  \lesssim \Lxpnorm{\Delta u}{3/2} \Lxpnorm{\nabla v^{+}}{2}^{2}
  \\
  &\lesssim \left( \sum\nolimits_{\pm} \Lxpnorm{v^{\pm}}{2}\Lxpnorm{v^{\pm}}{6} \right)
  \Lxpnorm{\nabla v^{+}}{2}^{2},
\end{align*}
where we used Lemma \ref{EllipticLemma}\eqref{Elliptic2} (see section 
\ref{Inequalities}) and the Sobolev embedding \eqref{L6Sobolev}.
Therefore, by Gronwall's lemma applied to $f(t) = \twonorm{\nabla 
v^{+}(t)}{}^{2}$,
\begin{equation}\label{vGradient}
  \twonorm{\nabla v^{+}(t)}{} \le \twonorm{\nabla v^{+}(0)}{} \exp\left( 
  \sum\nolimits_{\pm} \twonorm{v^{\pm}(0)}{}
  \int_{0}^{t} \Lpnorm{v^{\pm}(s)}{6} \, ds \right),
\end{equation}
where we used \eqref{L2conservation}. Therefore, \eqref{PSH1bound} will 
certainly follow if we can control the norms $\mixednormlocal{v^{\pm}}{2}{6}$.
To this end, define
$$
  Z_{T}^{\pm} = 
  \mixednormlocal{v^{\pm}}{\infty}{2}
  + \mixednormlocal{v^{\pm}}{2}{6}.
$$
In view of \eqref{L2conservation}, the second 
term can in fact be replaced by $\twonorm{v^{\pm}(0)}{}$. Then set $Z_{T} = 
Z_{T}^{+} + Z_{T}^{-}$.

We claim that (recall the notational conventions described in section 
\ref{MainResults})
\begin{equation}\label{ZineqforPS}
  Z_{T} \lesssim Z_{0}
  + \sigma(T) Z_{0}^{2} Z_{T}.
\end{equation}
This would imply that $Z_{T} \lesssim Z_{0}$ up to a time $T > 0$ 
only depending on $Z_{0} = \sum\nolimits_{\pm} \twonorm{v^{\pm}(0)}{}$. Then local 
well-posedness of \eqref{PS} in $L^{2}$ follows by standard arguments,\footnote{That is,
by exploiting the multilinearity of the nonlinear 
terms, the same argument gives estimates for a difference of two 
solutions in the norm $Z_{T}$. Then one can use, e.g., Picard iteration etc.}
hence global well-posedness by $L^{2}$-conservation.

So it remains to prove \eqref{ZineqforPS}.
To this end, we use a Strichartz type inequality for the
Schr\"odinger initial value problem on $\R^{1+3}$,
\begin{equation}\label{SchrIVP}
  i \partial_{t} v \pm \frac{\Delta}{2} v = F, \quad v \init = f.
\end{equation}
In fact, by Corollary 1.4 in \cite{KT}, if $2 \le q,r \le \infty$ and
$\tfrac{2}{q} + \tfrac{3}{r} = \tfrac{3}{2}$, then the 
estimate\footnote{This is the endpoint estimate, but one could also 
work with a non-endpoint norm $\mixed{2^{+}}{6^{-}}$ on the left hand 
side. This requires a modification of \eqref{vGradient}, of course.}
\begin{equation}\label{L2L6Str}
  \mixednormlocal{v}{2}{6} + \mixednormlocal{v}{\infty}{2}
  \lesssim \twonorm{f}{} + \mixednormlocal{F}{q'}{r'}
\end{equation}
holds for solutions of \eqref{SchrIVP}, where
$1 = \tfrac{1}{q} + \tfrac{1}{q'}$ and
$1 = \tfrac{1}{r} + \tfrac{1}{r'}$.
We apply this inequality with $q,r$ given by
$\tfrac{2}{q} = \varepsilon$ and $\tfrac{1}{r} = \tfrac{1}{2} - 
\tfrac{\varepsilon}{3}$, where $\varepsilon > 0$ is sufficiently small.
Thus, $(q',r') = (1^{+},2^{-})$. Applying \eqref{L2L6Str} to 
\eqref{PSb} then gives
$$
  Z_{T}^{\pm} \lesssim Z_{0}^{\pm} + \mixednormlocal{u}{q'}{3/\varepsilon}
  \Lxpnorm{v^{\pm}(0)}{2}.
$$
But using Sobolev embedding and \eqref{RieszCorollary}, followed by
 H\"older's inequality and $L^{p}$ interpolation,
$$
  \Lxpnorm{u}{3/\varepsilon}
  \lesssim \Lxpnorm{\Delta u}{(3/2)^{-}}
  \lesssim \sum\nolimits_{\pm} \Lxpnorm{v^{\pm}}{2}
  \Lxpnorm{v^{\pm}}{6^{-}}
  \lesssim \sum\nolimits_{\pm} \Lxpnorm{v^{\pm}}{2}^{1^{+}}
  \Lxpnorm{v^{\pm}}{6}^{1^{-}},
$$
and applying H\"older's inequality in $t$ then yields
$$
  \mixednormlocal{u}{q'}{3/\varepsilon} \lesssim \sigma(T)
  \sum\nolimits_{\pm} \Lxpnorm{v^{\pm}(0)}{2}^{1^{+}}
  \mixednormlocal{v^{\pm}}{2}{6}^{1^{-}} \lesssim \sigma(T) Z_{0} Z_{T}.
$$
This proves \eqref{ZineqforPS}.

\subsection{Some inequalities}\label{Inequalities}

Here we collect some simple estimates that will be used in later sections. 
First, for the operator $M$ defined by \eqref{Moperator}, we have:

\begin{lemma}\label{Mestimates}
The following operator norm estimates hold, for all $s \in \R$.
\begin{enumerate}
  \item\label{M2} $\norm{M^{-1}}_{H^{s} \to H^{s}} = O(1/c^{2})$.
  \item\label{M1} $\norm{M^{-1}}_{H^{s} \to H^{s+1}} = O(1/c)$.
  \item\label{M3} $\norm{M-c^{2}}_{H^{s+1} \to H^{s}} = O(c)$.
  \item\label{M4} $\norm{M-c^{2}}_{H^{s+2} \to H^{s}} = O(1)$. 
\end{enumerate}
\end{lemma}

\begin{proof} These statements translate to estimates on the 
Fourier symbols of the operators. Thus, the symbol $(c^{4} + c^{2} 
\abs{\xi}^{2})^{-1/2}$ of $M^{-1}$ is bounded by $c^{-2}$ as well as
$(c\abs{\xi})^{-1}$, which proves \eqref{M2} and \eqref{M1}, 
respectively.
The symbol of  $M - c^{2}$, given by \eqref{Mminusc2Symbol}, 
is bounded by $c\abs{\xi}$, and also by $\abs{\xi}^{2}/2$,
proving parts \eqref{M3} and \eqref{M4}, respectively.
\end{proof}

For the splitting \eqref{LowHigh} into low and high frequencies, we have:
\begin{lemma}\label{LowHighLemma}
The following estimates hold on $\R^{3}$.
\begin{enumerate}
  \item\label{HL1} $\Lpnorm{\tfrac{M}{c^{2}} f_{l}}{p} \lesssim 
  \Lpnorm{f_{l}}{p}$ for $1 \le p \le \infty$.
  \item\label{HL2} $\Sobnorm{f_{l}}{1+\varepsilon}
  \lesssim c^{\varepsilon} \Sobnorm{f_{l}}{1}$ for $\varepsilon > 0$.
  \item\label{HL3} $\Lpnorm{\tfrac{M}{c^{2}} f_{h}}{2} \lesssim 
  \tfrac{1}{c} \Sobdotnorm{f_{h}}{1}$.
  \item\label{HL4} $\Lpnorm{f_{h}}{2} \lesssim \tfrac{1}{c} 
  \Sobnorm{f_{h}}{1}$.
\end{enumerate}
\end{lemma}

\begin{proof}
Since $\tfrac{M}{c^{2}} f_{l} = \omega_{c} * f_{l}$, where 
$\widehat{\omega_{c}}(\xi) = \bigl(1 + 
\bigabs{\tfrac{\xi}{c}}^{2}\bigr)^{1/2} \chi\bigl(\tfrac{\xi}{2c}\bigr)$,
and since the $L^{1}$ norm of $\omega_{c}$ is independent of $c$, we get
\eqref{HL1} by Young's inequality.
The remaining inequalities are easy to prove using Plancherel's theorem; we omit 
the details.
\end{proof}

In order to estimate $A_{0}$, we will need:

\begin{lemma}\label{EllipticLemma}
The following estimates hold on $\R^{3}$.
\begin{enumerate}
  \item\label{Elliptic1} $\inftynorm{f}{}
  \lesssim \Lpnorm{\Delta f}{(3/2)^{+}}
  + \Lpnorm{\Delta f}{(3/2)^{-}}$.
  \item\label{Elliptic2} $\Lpnorm{\nabla f}{3}
  \lesssim \Lpnorm{\Delta f}{3/2}$.
\end{enumerate}
\end{lemma}

\begin{proof}
The second inequality is immediate from Sobolev embedding and the inequality
\eqref{RieszCorollary}.
To prove \eqref{Elliptic1}, observe first that for $\delta > 0$ 
arbitrarily small,
$$
  \inftynorm{f}{}
  \lesssim \Lpnorm{(I-\Delta)^{\delta} f}{3/\delta}
  \lesssim \Lpnorm{(-\Delta)^{\delta} f}{3/\delta}
  + \Lpnorm{f}{3/\delta}.
$$
The first inequality follows by Sobolev embedding (see \cite[Theorem 0.3.7]{So93}),
the second from \cite[Lemma V.2(ii)]{St}. By the Hardy-Littlewood-Sobolev
inequality (see \cite[Theorem 0.3.2]{So93}) the right hand side is
$\lesssim \Lpnorm{\Delta f}{3/(2-\delta)} + \Lpnorm{\Delta f}{3/(2+\delta)}$.
This concludes the proof.
\end{proof}

Finally, we note that the Sobolev embedding
\begin{equation}\label{L6Sobolev}
  \Lxpnorm{f}{6} \lesssim \Sobdotnorm{f}{1},
\end{equation}
implies
\begin{equation}\label{L2trilinear}
  \twonorm{fgh}{_{x}} \lesssim \Sobdotnorm{f}{1}
  \Sobdotnorm{g}{1} \Sobdotnorm{h}{1}
\end{equation}
and
\begin{equation}\label{L2bilinear}
  \twonorm{fg}{_{x}} \lesssim \Sobdotnorm{f}{1}
  \twonorm{g}{_x}^{1/2} 
  \Sobdotnorm{g}{1}^{1/2}.
\end{equation}
To prove the latter, write $\twonorm{fg}{}
\le \Lpnorm{f}{6} \Lpnorm{g}{3}$ and $\Lpnorm{g}{3} \le 
\twonorm{g}{}^{1/2} \Lpnorm{g}{6}^{1/2}$.

\section{Energy conservation and uniform $L^{2}$ 
bounds}\label{EnergyEstimates}

Throughout this section it is assumed that the hypotheses of Theorem \ref{Thm1} are 
satisfied. Our aim here is to prove the global-in-time $\mixed{\infty}{2}$ 
bound \eqref{L2bound} for $\phi^{\pm}$.
But by \eqref{PhiSplitting}, Lemma \ref{Mestimates}\eqref{M2} and
\eqref{L2bilinear},
$$
  \twonorm{\phi^{\pm}}{} \lesssim
  \twonorm{\phi}{}
  + \tfrac{1}{c^{2}} \twonorm{\partial_{t} \phi}{}
  + \tfrac{1}{c^{2}} \twonorm{\nabla A_{0}}{}
  \twonorm{\phi}{}^{1/2} \twonorm{\nabla \phi}{}^{1/2},
$$
at each time $t$, so it suffices to prove
\begin{equation}\label{EnergyBounds}
  \mixednorm{\phi}{\infty}{2}
  + \tfrac{1}{c} \mixednorm{\nabla \phi}{\infty}{2}
  + \tfrac{1}{c^{2}} \mixednorm{\partial_{t} \phi}{\infty}{2}
  + \tfrac{1}{c} \mixednorm{\nabla A_{0}}{\infty}{2} = O(1).
\end{equation}
This will be deduced from the conservation of the KGM energy $\mathcal 
E(t)$ given by \eqref{ConsLaw} and \eqref{T00}. Thus, if we can show
\begin{equation}\label{InitialEnergyBound}
  \mathcal E(0) = O(c^{2})
\end{equation}
and
\begin{equation}\label{EnergyControlsData}
  c^{2} \twonorm{\phi}{}^{2}
  + \twonorm{\nabla \phi}{}^{2}
  + \tfrac{1}{c^{2}} \twonorm{\partial_{t} \phi}{}^{2}
  + \twonorm{\nabla A_{0}}{}^{2}
  \lesssim \mathcal E\left( 1 + \mathcal E/c^3 + \mathcal E^2 / c^6 \right)
\end{equation}
at each time $t$, then \eqref{EnergyBounds} follows immediately.

\subsection{Proof of \eqref{InitialEnergyBound}}
In view of the definitions \eqref{ConsLaw}, \eqref{T00}, \eqref{Ds} and
\eqref{EandB}, it is enough to prove, at $t = 0$,
\begin{align}
  \label{InitialEnergyBoundA}
  c^{2} \twonorm{\phi}{}^{2}
  + \twonorm{\nabla \phi}{}^{2}
  + \tfrac{1}{c^{2}} \twonorm{\partial_{t} \phi}{}^{2} &= O(c^{2})
  \\
  \label{InitialEnergyBoundC}
  \twonorm{\nabla A_{0}}{}^{2} + \twonorm{\nabla \A}{}^{2}
  + \tfrac{1}{c^{2}} \twonorm{\partial_{t} \A}{}^{2} &= O(c^{2}),
  \\
  \label{InitialEnergyBoundB}
  \tfrac{1}{c^{2}} \twonorm{A_{0} \phi}{}^{2}
  + \tfrac{1}{c^{2}} \twonorm{\A \phi}{}^{2} &= O(1).
\end{align}
The first two terms on LHS\eqref{InitialEnergyBoundA} are $O(c^{2})$ 
at $t = 0$
by \eqref{PhidataBound}, and for the third term we write $\partial_{t} \phi
= M M^{-1} \partial_{t} \phi$, which gives
\begin{equation}\label{TimeDerBound}
  \tfrac{1}{c} \twonorm{\partial_{t} \phi}{} \lesssim c 
  \twonorm{M^{-1}\partial_{t} \phi}{} + \Sobnorm{M^{-1} \partial_{t} 
  \phi}{1} = O(c)
\end{equation}
by \eqref{PhidataBound}. This proves \eqref{InitialEnergyBoundA}.

The last two terms on LHS\eqref{InitialEnergyBoundC} are $O(1)$ at 
$t = 0$
by \eqref{AdataBound}, and for the first term we use the elliptic 
estimate (see \cite[Eqs. (3.4a,b)]{KM}) $\twonorm{\nabla A_{0}(t)}{} 
\lesssim \tfrac{1}{c} \twonorm{\partial_{t} \phi}{}$. Therefore, by 
\eqref{TimeDerBound}, $\twonorm{\nabla A_{0}(t=0)}{} = O(c)$, and 
this concludes the proof of \eqref{InitialEnergyBoundC}.

Finally, to prove \eqref{InitialEnergyBoundB} at $t = 0$, use \eqref{L2bilinear}
and the bounds in \eqref{PhidataBound} and 
\eqref{InitialEnergyBoundC}.

\subsection{Proof of \eqref{EnergyControlsData}}
First, by \cite[Eq. (1.3c)]{KM},
\begin{equation}\label{EnergyControlsDataA}
  \twonorm{\nabla A_{0}}{}^{2} + \twonorm{\nabla \A}{}^{2} \lesssim 
  \mathcal E
\end{equation}
for all $t$, so we get the desired bound for the last term on 
LHS\eqref{EnergyControlsData}. The first term is obviously bounded 
by $\mathcal E$, so it remains to consider the two middle terms.
But using the definition \eqref{Ds} and \eqref{L2bilinear},
$$
  \twonorm{\nabla \phi}{} \le \sum\nolimits_{1}^{3} \twonorm{D_{j} \phi}{}
  + \tfrac{1}{c} \twonorm{\nabla \A}{} \twonorm{\phi}{}^{1/2} 
  \twonorm{\nabla \phi}{}^{1/2}.
$$
Now use the fact that if $\alpha \le \beta + \gamma \sqrt \alpha$, 
where $\alpha, \beta, \gamma \ge 0$, then $\alpha \le 2\beta + 4 
\gamma^{2}$. Combining this with \eqref{EnergyControlsDataA} gives the 
bound $\twonorm{\nabla \phi}{}^{2} \lesssim \mathcal E + (\mathcal 
E/c^{2})^{3}$. Similarly,
$$
  \tfrac{1}{c} \twonorm{\partial_{t} \phi}{} \le \twonorm{D_{0} \phi}{}
  + \tfrac{1}{c} \twonorm{\nabla A_{0}}{} \twonorm{\phi}{}^{1/2} 
  \twonorm{\nabla \phi}{}^{1/2}.
$$
Squaring this, and using \eqref{EnergyControlsDataA} as well as the 
bounds already obtained for $\phi$ and $\nabla \phi$, we get the 
correct bound for $\tfrac{1}{c^{2}} \twonorm{\partial_{t} \phi}{}^{2}$.

\section{Linear and bilinear spacetime 
estimates}\label{SpacetimeEstimates}

Here we prove some linear and bilinear Strichartz type estimates on 
$\R^{1+3}$ for the operators $L^{\pm}$, defined by \eqref{Loperator}.

\subsection{Linear estimates}
The key observation is that the propagators associated to $L^{\pm}$,
\begin{equation}\label{Lpropagator}
  U^{\pm}(t) = e^{\mp it(M-c^{2})},
\end{equation}
behave like the Schr\"odinger propagators
\begin{equation}\label{SchrodingerPropagator}
  V^{\pm}(t) = e^{\pm it\Delta/2}
\end{equation}
at low frequencies ($\lesssim c$) and like the wave equation propagators
$e^{\mp itc\sqrt{-\Delta}}$ at high frequencies ($\gg c$). Indeed,
$U^{\pm}(t)$ is a multiplier with Fourier symbol
$e^{\mp it h_{c}(\xi)}$, where
\begin{equation}\label{hc}
  h_{c}(\xi) = \frac{\abs{\xi}^{2}}{1 + \sqrt{1 + \abs{\xi}^{2}/c^{2}}}
  \sim
  \begin{cases}
    \abs{\xi}^{2}/2 &\quad \text{for $\abs{\xi} \lesssim c$},
    \\
    c\abs{\xi} &\quad \text{for $\abs{\xi} \gg c$}.
  \end{cases}
\end{equation}
It is therefore not surprising that we have Strichartz estimates for 
$U^{\pm}$ in $\mixed{q}{r}$ for every sharp wave admissible pair
$(q,r)$ of Lebesgue exponents, and if we restrict to low frequency
($\lesssim c$), Schr\"odinger admissible exponents are also 
allowed.

Let us be more explicit. Following the terminology introduced in 
\cite{KT}, we say that a pair $(q,r)$ of Lebesgue exponents is \emph{sharp 
wave admissible} (for $\R^{1+3}$) if
\begin{equation}\label{WaveAdm}
  \frac{1}{q} + \frac{1}{r} = \frac{1}{2} \quad \text{and} \quad
  (q,r) \neq (2,\infty),
\end{equation}
and we say $(q,r)$ is \emph{Schr\"odinger admissible} (for 
$\R^{1+3}$) if $q,r \ge 2$ and
\begin{equation}\label{SchrodingerAdm}
  \frac{2}{q} + \frac{3}{r} = \frac{3}{2}.
\end{equation}

\begin{proposition}\label{WaveProp}
For every sharp wave admissible pair $(q,r)$, the estimate
\begin{equation}\label{WaveEst}
  \mixednormlocal{U^{\pm}(t) f}{q}{r}
  \lesssim \Sobdotnorm{f}{\frac{1}{q}}
  + c^{-\frac{1}{q}} \Sobdotnorm{f}{\frac{2}{q}}
\end{equation}
holds.
\end{proposition}

The choice of norm on the right hand side is motivated by dimensional
analysis. Thus, the first term $\Sobdotnorm{f}{\frac{1}{q}}$, which 
dominates at low frequency, is what one would get by scaling if $U^{\pm}$ 
were replaced by the Schr\"odinger propagator $V^{\pm}$. If instead we 
consider high frequencies and replace $U^{\pm}$ by the wave 
propagator $e^{itc\sqrt{-\Delta}}$, we get the second term
$c^{-\frac{1}{q}} \Sobdotnorm{f}{\frac{2}{q}}$, again by scaling.

Then using Duhamel's principle to write the solution of
\begin{equation}\label{IVPforL}
  L^{\pm} u = F, \quad u \init = f
\end{equation}
as
\begin{equation}\label{Lsolution}
  u(t) = U^{\pm}(t) f + \int_{0}^{t} U^{\pm}(t-s) F(s) \, ds,
\end{equation}
and noting that the norm on RHS\eqref{WaveEst} is dominated by
$\Sobnorm{f}{\frac{2}{q}}$ as $c \to \infty$, we immediately obtain the following:

\begin{corollary}
For every sharp wave admissible pair $(q,r)$, the estimate
\begin{equation}\label{WaveEstPrime}
  \mixednormlocal{u}{q}{r} + \LpHslocal{u}{\infty}{\frac{2}{q}} 
  \lesssim \Sobnorm{f}{\frac{2}{q}}
  + \int_{0}^{T} \Sobnorm{F(t)}{\frac{2}{q}} \, dt
\end{equation}
holds for solutions of \eqref{IVPforL}.
\end{corollary}

Next we consider estimates for Schr\"odinger admissible exponents.

\begin{proposition}\label{SchrodingerProp}
Let $(q,r)$ and $(\widetilde q, \widetilde r)$ be any two 
Schr\"odinger admissible pairs. Then for the low frequency part $u_{l}$
(see \eqref{LowHigh} for definition) of the solution of  \eqref{IVPforL} we have the estimate
$$
  \mixednormlocal{u_{l}}{q}{r} + \mixednormlocal{u_{l}}{\infty}{2}
  \lesssim \twonorm{f_{l}}{} + \mixednormlocal{F_{l}}{\widetilde 
  q'}{\widetilde r'},
$$
where $\tfrac{1}{\widetilde q} + \tfrac{1}{\widetilde q'} = 1$ and
$\tfrac{1}{\widetilde r} + \tfrac{1}{\widetilde r'} = 1$.
\end{proposition}

Let us turn to the proofs.

\subsubsection{Proof of Proposition \ref{WaveProp}}

Proceeding as in the standard proof of the Strichartz estimates for
the homogeneous wave equation (see, e.g., \cite{KT} or
\cite[Section III.5]{So}) we reduce to proving the decay estimate
\begin{equation}\label{KmuEst}
  \abs{K_{\mu,c}(t,x)} \lesssim
  \begin{cases}
    \frac{\mu}{\abs{t}} &\quad \text{for $\mu \lesssim c$},
    \\
    \frac{\mu^{2}}{c\abs{t}} &\quad \text{for $\mu \gg c$},
  \end{cases}
\end{equation}
for the convolution kernel
$$
  K_{\mu,c}(t,x) = \int_{\R^{3}} e^{ix\cdot \xi} 
  e^{ith_{c}(\xi)} \beta\bigl(\tfrac{\xi}{\mu}\bigr) \, d\xi,
$$
where $h_{c}$ is given by \eqref{hc}, $\beta$ is a 
Littlewood-Paley cut-off function supported in the annulus 
$\abs{\xi} \sim 1$ and $\mu$ is a dyadic number of the form $2^{j}$, $j 
\in \mathbb Z$. But in view of the scaling identity
$$
  K_{\mu,c} (t,x) = c^{3} K_{\mu/c,1}(c^{2}t,cx),
$$
it suffices to prove \eqref{KmuEst} for $c = 1$. To simplify the 
notation we write $K_{\mu} = K_{\mu,1}$ and $h = h_{1}$. We shall need 
the following fact, whose elementary proof we omit:

\begin{lemma}\label{alphaLemma}
  Define $\alpha(r) = \frac{r^{2}}{ 1 + \sqrt{1+ r^{2}}}$
  for $r > 0$. Then
  $\alpha'(r) = \frac{r}{\sqrt{1+r^{2}}}$ and
  $\alpha''(r) = \frac{1}{(1+r^{2})^{3/2}}$.
\end{lemma}

To prove \eqref{KmuEst} for $c = 1$, we split into four cases:
\begin{enumerate}
  \item\label{case1}
  $\mu \lesssim 1$ and $\abs{x} \gtrsim \mu \abs{t}$, 
  \item\label{case2}
  $\mu \lesssim 1$ and $\abs{x} \ll \mu \abs{t}$, 
  \item\label{case3}
  $\mu \gg 1$ and $\abs{x} \gtrsim \abs{t}$, 
  \item\label{case4}
  $\mu \gg 1$ and $\abs{x} \ll \abs{t}$.
\end{enumerate}
Introducing polar coordinates $\xi = r\omega$, $r > 0$, $\omega \in 
S^{2}$, we have
\begin{align}
  \label{Polar1}
  K_{\mu}(t,x) &= \int_{0}^{\infty} \int_{S^{2}}  e^{i r x \cdot \omega} e^{it 
  \alpha(r)} \beta\bigl(\tfrac{r}{\mu}\bigr) r^{2} \, d\sigma(\omega) \, dr
  \\
  \label{Polar2}
  &= \int_{0}^{\infty} \widehat \sigma(rx) e^{it 
  \alpha(r)} \beta\bigl(\tfrac{r}{\mu}\bigr) r^{2} \, dr,
\end{align}
where $\sigma$ is surface measure on $S^{2}$. Since $\abs{\widehat 
\sigma (\xi)} \lesssim \abs{\xi}^{-1}$ (see, e.g., \cite[Eq. (5.13)]{So})
we get from \eqref{Polar2}
$$
  \abs{K_{\mu}(t,x)} \lesssim \abs{x}^{-1} \int_{0}^{\infty} 
  \beta\bigl(\tfrac{r}{\mu}\bigr) r \, dr \sim \frac{\mu^{2}}{\abs{x}},
$$
which proves \eqref{KmuEst} ($c=1$) for the cases \eqref{case1} and 
\eqref{case3}. Next, rewrite \eqref{Polar1} as
$K_{\mu}(t,x) = \int_{S^{2}} I(\omega) \, d\sigma(\omega)$, where
$$
  I(\omega) = \int_{0}^{\infty} \frac{d}{dr} \left[ e^{i (t 
  \alpha(r) + r x\cdot \omega)} \right]
  \frac{\beta\bigl(\tfrac{r}{\mu}\bigr) r^{2}}{i \bigl( t \alpha'(r) + 
  x \cdot \omega \bigr)} \, dr.
$$
Integrating by parts and writing
$$
  - \frac{d}{dr} \left[ \frac{\beta\bigl(\tfrac{r}{\mu}\bigr) r^{2}}
  {i \bigl( t \alpha'(r) +   x \cdot \omega \bigr)} \right]
  =  \frac{\beta\bigl(\tfrac{r}{\mu}\bigr) r^{2} t \alpha''(r)}
  {i \bigl( t \alpha'(r) + x \cdot \omega \bigr)^{2}}
  - \frac{\tfrac{d}{dr} \bigl[ \beta\bigl(\tfrac{r}{\mu}\bigr) r^{2} 
  \bigr]}{ i \bigl( t \alpha'(r) + 
  x \cdot \omega \bigr)}
$$
gives $I = I_{1} + I_{2}$.

Consider case \eqref{case2}. Then $r \sim \mu \lesssim 1$, so
$\abs{\alpha'(r)} \sim \mu$ and $\abs{\alpha''(r)} \sim 1$
by Lemma \ref{alphaLemma}.
Then since $\abs{x} \ll \mu \abs{t}$, we get $\abs{t \alpha'(r) + 
x \cdot \omega} \gtrsim \mu \abs{t}$, and this gives
$\abs{I_{j}(\omega)} \lesssim \mu/\abs{t}$ for $j = 1,2$, proving
\eqref{KmuEst} for this case.

Finally, consider case \eqref{case4}. Then Lemma \ref{alphaLemma}
gives $\abs{\alpha'(r)} \sim 1$ and $\abs{\alpha''(r)}
\sim \mu^{-3}$, since $r \sim \mu \gg 1$. In view of the assumption
$\abs{x} \ll \abs{t}$, we then get
$\abs{t \alpha'(r) + x \cdot \omega} \gtrsim \abs{t}$, whence
$\abs{I_{j}(\omega)} \lesssim \mu^{2}/\abs{t}$ for $j = 1,2$.
This proves \eqref{KmuEst} for case \eqref{case4}, and concludes the
proof of Proposition \ref{WaveProp}.

\subsubsection{Proof of Proposition \ref{SchrodingerProp}}

Take the convolution with $\theta_{c}$ in \eqref{IVPforL} and use 
the identity $\theta_{c} = \theta_{c} * \theta_{2c}$ to see that
$L^{\pm}_{\text{low}} u_{l} = F_{l}$ with data $u_{l} \init = f_{l}$, 
where $L^{\pm}_{\text{low}}$ is the operator with propagator
$U^{\pm}_{\text{low}}(t) = \theta_{2c} * e^{\mp it (M-c^{2})}$.
It therefore suffices to prove
$$
  \mixednormlocal{u}{q}{r} + \mixednormlocal{u}{\infty}{2}
  \lesssim \twonorm{f}{} + \mixednormlocal{F}{\widetilde 
  q'}{\widetilde r'}
$$
for solutions of $L^{\pm}_{\text{low}} u = F$ with data $u \init = f$.
But by \cite[Theorem 1.2]{KT} (see also the proof of Corollary 1.4 there)
it suffices to prove the decay estimate
\begin{equation}\label{Kest}
  \abs{K_c(t,x)} \lesssim \abs{t}^{-3/2}
\end{equation}
for the convolution kernel $K_c(t,x) = (2\pi)^{-3} \int_{\R^{3}} 
e^{i x\cdot \xi} e^{ith_{c}(\xi)} \chi(\xi/c) \, d\xi$, where $h_{c}$ 
is given by \eqref{hc}. In view of the scaling identity $K_c(t,x) = c^{3} 
K_1(c^{2}t,cx)$, it is enough to prove \eqref{Kest} for $c = 1$, in which case it
follows from a standard result about decay of the Fourier transform of surface 
carried measures; see \cite[Theorem 1.2.1]{So93}. Indeed, $K_1(t,x)$
is the (inverse) spacetime Fourier transform of the measure
(recall that $h = h_{1}$)
$$
  \delta( \tau - h(\xi) ) \chi(\xi),
$$
which is compactly supported on the hypersurface $\{ (\tau,\xi) \in 
\R^{1+3} : \tau = h(\xi) \}$, whose curvature is non-vanishing.

\subsection{Bilinear null form estimates}

In \cite{KM93}, Klainerman and Machedon proved that the estimate
\begin{equation}\label{False}
  \twonorm{u \nabla v}{(\R^{1+3})} \lesssim \Sobnorm{f}{1} \Sobnorm{g}{1}
\end{equation}
fails for solutions of $\square u = \square v = 0$ on $\R^{1+3}$ with initial data 
$(f,0)$ and $(g,0)$. In particular, this shows that the endpoint 
$(q,r) = (2,\infty)$ for the linear Strichartz estimates is forbidden,
for if the estimate $\mixednorm{u}{2}{\infty} \lesssim 
\Sobnorm{f}{1}$ were true, it would clearly imply \eqref{False}.
If the bilinear form $u \nabla v$ in \eqref{False} is replaced by one of the null 
forms $Q_{ij}\bigl(\D^{-1} u,v\bigr)$ or
$\D^{-1} Q_{ij}(u,v)$, the estimate is true, however, as proved in \cite{KM93}.
Here $\D^\alpha = (-\Delta)^{\alpha/2}$ and
$$
  Q_{ij}(u,v) = \partial_{i} u \, \partial_{j} v - \partial_{j} u \, \partial_{i} v,
  \quad (1 \le i,j \le 3).
$$
This fact was used in \cite{KM} to control the bilinear terms with 
derivatives in the KGM system, which turn out to have this 
structure when the Coulomb gauge \eqref{CG} condition is used.

In fact (see the proof of the corollary to Proposition 2.1 in \cite{KM})
\begin{equation}\label{EasyNullFormControl}
  \twonorm{\Proj(u \nabla v)}{_{x}} \lesssim \sum_{1 \le i,j \le 3}
  \twonorm{\D^{-1} Q_{ij}(u,v)}{_{x}},
\end{equation}
where the projection $\Proj$ is given by \eqref{ProjDef}. Moreover, if 
$u = (u^{1}, u^{2}, u^{3})$ is vector valued and
divergence free, so that $\Proj u = u$, then (see the proof of Proposition 2.2
in \cite{KM}, or \cite[Section 1.5]{Se})
\begin{equation}\label{HardNullFormIdentity}
  u \cdot \nabla v = \half \sum_{1 \le i,j \le 3} Q_{ij}\bigl(\D^{-1} [ 
  R_{j} u^{i} - R_{i} u^{j}], v \bigr),
\end{equation}
where $R_{i} = \D^{-1} \partial_{i}$ are the Riesz operators.

Here we prove versions of the Klainerman-Machedon null form estimates where one
or both of $u,v$ solve $L^{\pm} [\cdot] = 0$ instead of $\square [\cdot] = 0$.

\begin{proposition}\label{EasyNullFormProp}
Suppose $L^{\pm} u = F$ and $L^{\pm} v = G$ (independent signs) with 
initial data $u \init = f$ and $v \init = g$. Then
$$
  \twonorm{\D^{-1} Q(u,v)}{(S_{T})}
  \lesssim \left( \Sobnorm{f}{1} + \LpHslocal{F}{1}{1} \right)
  \left( \Sobnorm{g}{1} + \LpHslocal{G}{1}{1} \right)
$$
for $Q = Q_{ij}$.
\end{proposition}

In view of \eqref{EasyNullFormControl}, this implies the following:

\begin{corollary} Under the hypotheses of Proposition 
\ref{EasyNullFormProp}, we have
$$
  \twonorm{\Proj(u \nabla v)}{(S_{T})}
  \lesssim \left( \Sobnorm{f}{1} + \LpHslocal{F}{1}{1} \right)
  \left( \Sobnorm{g}{1} + \LpHslocal{G}{1}{1} \right).
$$
\end{corollary}

Next, we consider the null form $Q_{ij}\bigl(\D^{-1} u,v\bigr)$.

\begin{proposition}\label{HardNullFormProp}
Suppose $\square u = F$ and $L^{\pm} v = G$ with initial data
$u \init = f_{0}$, $\partial_{t} u \init = f_{1}$ and $v \init = g$. Then
\begin{multline*}
  \twonorm{Q\bigl(\D^{-1} u,v\bigr)}{(S_{T})}
  \\
  \lesssim c^{-1/2}\left( \Sobdotnorm{f_{0}}{1} + \tfrac{1}{c} 
  \twonorm{f_{1}}{} + c \mixednormlocal{F}{1}{2} \right)
  \left( \Sobdotnorm{g}{1} + \norm{G}_{L_{t}^{1}\dot H^{1}(S_{T})} \right)
\end{multline*}
for $Q = Q_{ij}$.
\end{proposition}

Then using \eqref{HardNullFormIdentity} and noting that the Riesz 
operators commute with $\square$ and are bounded on every $H^{s}$ 
space, we obtain:

\begin{corollary} Assume the hypotheses of Proposition 
\ref{HardNullFormProp} are satisfied. If in addition we 
assume that $u(t)$ is vector valued and divergence free, then
\begin{multline*}
  \twonorm{u \cdot \nabla v}{(S_{T})}
  \\
  \lesssim c^{-1/2}\left( \Sobdotnorm{f_{0}}{1} + \tfrac{1}{c} 
  \twonorm{f_{1}}{} + c \mixednormlocal{F}{1}{2} \right)
  \left( \Sobdotnorm{g}{1} + \norm{G}_{L_{t}^{1}\dot H^{1}(S_{T})} \right).
\end{multline*}
\end{corollary}

In the rest of this section, the Fourier transform of a function 
$u(t,x)$ (resp. $f(x)$) is denoted $\Fourier u(\tau,\xi)$ (resp. 
$\widehat f(\xi)$). Then
\begin{equation}\label{QijFT}
  \Fourier [ Q_{ij}(u,v) ](\tau,\xi)
  = \int q_{ij}(\eta,\xi-\eta) 
  \Fourier u(\lambda,\eta) \Fourier v(\tau-\lambda,\xi-\eta)\, d\lambda \, d\eta,
\end{equation}
where $q_{ij}(\xi,\eta) = \xi_{i} \eta_{j} - \xi_{j} \eta_{i}$ for 
$\xi, \eta \in \R^{3}$.
We will need the two inequalities
\begin{align}
  \label{QsymbolEst}
  \abs{q_{ij}(\xi,\eta)} \le \abs{\xi \times \eta} \le
  \abs{\xi + \eta} \abs{\xi}^{1/2} \abs{\eta}^{1/2}.
\end{align}
The first inequality is obvious, and to prove the second, observe 
that $\xi \times \eta = (\xi + \eta) \times \eta = \xi \times (\xi + 
\eta)$, whence $\abs{\xi \times \eta} \le \abs{\xi + \eta} \min 
(\abs{\xi}, \abs{\eta})$. From \eqref{QijFT}, \eqref{QsymbolEst} and 
Plancherel's theorem, we then get
\begin{equation}\label{EasyNullFormReduction}
  \twonorm{\D^{-1} Q_{ij}(u,v)}{(\R^{1+3})}
  \le \twonorm{\D^{1/2} u \cdot\D^{1/2} v}{(\R^{1+3})}  
\end{equation}
provided $\Fourier u, \Fourier v \ge 0$.

\subsection{Proof of Proposition \ref{EasyNullFormProp}}

In view of the formula \eqref{Lsolution} for the solution of 
\eqref{IVPforL}, it suffices to prove this for $F = G = 0$.
Then $\Fourier u (\tau,\xi) = \delta(\tau \pm h_{c}(\xi)) \widehat 
f(\xi)$ and $\Fourier v (\tau,\xi) = \delta(\tau \pm h_{c}(\xi)) \widehat 
g(\xi)$, where $h_{c}$ is given by \eqref{hc}.
Withouth loss of generality, we assume $\widehat f, \widehat g \ge 0$.
Thus, \eqref{EasyNullFormReduction} applies, and
$$
  \twonorm{\D^{-1} Q_{ij}(u,v)}{(\R^{1+3})}
  \le \norm{\D^{1/2} u}_{L^{4}(\R^{1+3})} \norm{\D^{1/2} 
  v}_{L^{4}(\R^{1+3})}.
$$
Now apply Proposition \ref{WaveProp}, or its corollary, with 
$q = r = 4$, to conclude the proof.

\subsection{Proof of Proposition \ref{HardNullFormProp}}

Reasoning as above, we may assume $G = 0$, so that $v(t) = 
U^{\pm}(t) g$. Similarly, since the solution of
\begin{equation}\label{IVPforWave}
  \square u = F, \quad u \init = f_{0}, \quad \partial_{t} u \init = f_{1}
\end{equation}
is given by (recall that $\square = - \tfrac{1}{c^{2}} 
\partial_{t}^{2} + \Delta$)
\begin{multline}\label{WaveSolution}
  u(t) = \cos(ct\D) f_{0} + (c\D)^{-1} \sin(ct\D) f_{1}
  \\
  - c \int_{0}^{t} \D^{-1} \sin\bigl(c(t-s)\D\bigr) F(s) \, ds,
\end{multline}
we reduce to the case where $f_{1} = 0$, $F = 0$ and $u(t) = e^{\pm 
ict\D} f_{0}$. Without loss of generality, we choose the plus sign in 
the exponential. Thus, writing $f = f_{0}$, we only have to prove
$$
  \twonorm{Q_{ij}\bigl(\D^{-1} u,v\bigr)}{(\R^{1+3})}
  \lesssim c^{-1/2} \Sobdotnorm{f}{1} \Sobdotnorm{g}{1},
$$
where $u(t) = e^{ict\D} f$ and $v(t) = e^{\pm it (M-c^{2})} g$.
Changing variables $t \to ct$, this becomes
\begin{equation}\label{ReducedEst}
  \twonorm{Q_{ij}\bigl(\D^{-1} u',v'\bigr)}{(\R^{1+3})}
  \lesssim \Sobdotnorm{f}{1} \Sobdotnorm{g}{1},
\end{equation}
where $u'(t) = e^{it\D} f$ and $v'(t) = e^{\pm it (M-c^{2})/c} g$.
Thus,
$$\Fourier u'(\tau,\xi) = \delta(\tau - \abs{\xi}) \widehat 
f(\xi) \quad \text{and}
\quad \Fourier v'(\tau,\xi) =
\delta\left(\tau \pm c \, h\bigl(\tfrac{\xi}{c}\bigr) \right)
\widehat g(\xi),
$$
where $h(\xi) = \tfrac{\abs{\xi}^{2}}{1 + \sqrt{1 + \abs{\xi}^{2}}}$.
We may assume $\widehat f, \widehat g \ge 0$. Then by \eqref{QijFT},
\begin{multline*}
  \abs{\Fourier [ Q_{ij}(\D^{-1}u',v') ](\tau,\xi)}
  \\
  \le \int \frac{\abs{q_{ij}(\eta,\xi-\eta)}}{\abs{\eta}^{2} \abs{\xi-\eta}}
  (\D f)\, \widehat \,(\eta)
  (\D g)\, \widehat \,(\xi-\eta)
  \delta\left(\tau - \abs{\eta} \pm c \, h\bigl(\tfrac{\xi-\eta}{c}\bigr) \right)
  \, d\eta.
\end{multline*}
Now apply the Cauchy-Schwarz inequality 
with respect to the measure $\delta(\dots) \, d\eta$, square both sides
and integrate in $d\tau \, d\xi$ to obtain
$$
  \twonorm{Q_{ij}\bigl(\D^{-1} u',v'\bigr)}{(\R^{1+3})}^{2}
  \le \inftynorm{I^{\pm}}{} \Sobdotnorm{f}{1}^{2} 
  \Sobdotnorm{g}{1}^{2},
$$
where
$$
  I^{\pm}(\tau,\xi) = \int \frac{\abs{q_{ij}(\eta,\xi-\eta)}^{2}}
  {\abs{\eta}^{4} \abs{\xi-\eta}^{2}}
  \delta\left(\tau - \abs{\eta} \pm c \, h\bigl(\tfrac{\xi-\eta}{c}\bigr) \right)
  \, d\eta.
$$
This reduces \eqref{ReducedEst} to proving that $I^{\pm}$ is 
bounded, independently of $c$. But by \eqref{QsymbolEst},
$$
  I^{\pm} \le \int \frac{\sin^2 \theta}
  {\abs{\eta}^{2}}
  \delta\left(\tau - \abs{\eta} \pm c \, h\bigl(\tfrac{\xi-\eta}{c}\bigr) \right)
  \, d\eta,
$$
where $\theta$ denotes the angle between $\eta$ and $\xi-\eta$.
Now apply the following general result, with $k(r) = c \alpha(r/c)$ 
and $\alpha$ as in Lemma \ref{alphaLemma}.

\begin{lemma}\label{DeltaLemma} Suppose $k(r)$ is positive and 
differentiable for $r > 0$, and that $\abs{k'(r)} \le 1$. Define
$$
  I^{\pm}(\tau,\xi) =  \int \frac{\sin^2 \theta}{\abs{\eta}^{2}}
  \delta\bigl(\tau - \abs{\eta} \pm k(\abs{\xi-\eta}) \bigr)
  \, d\eta,
$$
where $\theta$ is the angle between $\eta$ and $\xi-\eta$.
Then $\sup_{\tau, \xi} I^{\pm}(\tau,\xi) \le 8\pi$.
\end{lemma}

To see that this applies with $k(r) = c \alpha(r/c)$, we need only
observe that $k'(r) = \alpha'(r/c)$, and $0 < \alpha' < 1$ by Lemma 
\ref{alphaLemma}. We remark that the above lemma also applies with $k(r) = 
r$, which corresponds to the Klainerman-Machedon estimates (then $u$ and 
$v$ both solve the homogeneous wave equation).

\begin{proof}[Proof of Lemma \ref{DeltaLemma}]
In polar coordinates, $I^{\pm}(\tau,\xi) = \int_{S^{2}} \rho(\tau,\xi,\omega)
\, d\omega$, where
$$
  \rho(\tau,\xi,\omega) = \int_{0}^{\infty} (\sin \theta)^2
  \delta\bigl(\tau - r \pm k(\abs{\xi-r\omega}) \bigr) \, dr,
$$
so it suffices to show that $\rho \le 2$ for all $\tau, \xi$, and 
for almost every $\omega \in S^{2}$.

We shall use the following fact: Suppose $f : \R \to \R$ is differentiable
with $f'(r) < 0$, and has a zero at $r_{0}$. Then (see \cite[Theorem 6.1.5]{Horm})
\begin{equation}\label{DeltaIdentity}
  \delta\bigl( f(r) \bigr) \, dr = \frac{\delta(r-r_{0}) \, dr}{\abs{f'(r_{0})}}.
\end{equation}
Take $f(r) = \tau - r \pm k(\abs{\xi-r\omega})$, for 
fixed $\tau,\xi$ and $\omega$. Then, since $\abs{k'} \le 1$,
\begin{equation}\label{fDer}
  f'(r) = - 1 \pm k'(\abs{\xi-r\omega}) 
  \frac{(\xi-r\omega)\cdot\omega}{\abs{\xi-r\omega}}
  \le -1 + \abs{ \cos \theta},
\end{equation}
so $f'(r) < 0$ if we exclude the two points on $S^{2}$ where 
$\omega$ is collinear with $\xi$. Since \eqref{fDer} 
shows that $\abs{f'} \ge 1 - \abs{\cos \theta} \ge \tfrac{1}{2}\sin^{2} \theta$,
we conclude from \eqref{DeltaIdentity} that $\rho(\tau,\xi,\omega) \le 2$.
\end{proof}

\section{Local-in-time a priori bounds}\label{AprioriBounds}

Here we prove parts \eqref{Thm1A}--\eqref{Thm1H} of
Theorem \ref{Thm1}. Throughout this section we assume that the
hypotheses of the theorem are satisfied.

Let us first prove Lemma \ref{SplitSystemLemma}.
Solving \eqref{PhiSplitting} for $i\partial_{t} \phi$ gives
\begin{equation}\label{PhiTimeDer}
  i \partial_{t} \phi = M(\phi^{+} - \phi^{-}) + A_{0} (\phi^{+} + 
  \phi^{-}).
\end{equation}
Inserting this into \eqref{KGMa} gives \eqref{KGM'a}. Since 
\eqref{KGM'b} is the same as \eqref{KGMb}, it only remains to check 
\eqref{KGM'c}. To do this, take a time derivative of 
\eqref{PhiSplitting}, and use \eqref{KGMc} to eliminate 
$\partial_{t}^{2} \phi$. This gives
$$
  \partial_{t} \phi^{\pm} = \half \left\{ \partial_{t} \phi \mp iM 
  \phi \pm M^{-1} \left[ - 2c \A \cdot \nabla \phi + A_{0} 
  \partial_{t} \phi + i (A_{0}^{2} - \A^{2}) \phi \right] \right\}.
$$
Using \eqref{PhiTimeDer} to eliminate $\partial_{t} \phi$, the right 
hand side becomes
\begin{equation}\label{RHSversion1}
  \mp iM \phi^{\pm}
  - \frac{i}{2} A_{0} \phi
  \pm \half M^{-1} \left\{ - 2c \A \cdot \nabla \phi
  - i A_{0} M (\phi^{+} - \phi^{-})
  -i \A^{2} \phi \right\},
\end{equation}
and this can be rewritten as
\begin{equation}\label{RHSversion2}
  \mp iM \phi^{\pm}
  - i A_{0} \phi^{\pm}
  \pm \half M^{-1} \left\{ - 2c \A \cdot \nabla \phi
  - i [A_{0},M](\phi^{+} - \phi^{-})
  -i \A^{2} \phi \right\}.
\end{equation}
Adding $\pm i c^{2} \phi^{\pm}$ and then multiplying by $e^{\pm i tc^{2}}$
gives \eqref{KGM'c}, since $[A_{0},M] = [A_{0},M-c^{2}]$ and, by \eqref{PsiDef},
$$
  \partial_{t} \psi^{\pm} = e^{\pm i tc^{2}} \left( \partial_{t} 
  \phi^{\pm} \pm i c^{2} \phi^{\pm} \right).
$$

In fact, the commutator structure is not 
needed in Theorem \ref{Thm1} except to prove the bound for $R$.
To simplify certain arguments we will therefore
use an alternative formulation of \eqref{KGM'c}, obtained by using
the expression \eqref{RHSversion1} for
$\partial_{t} \phi^{\pm}$, instead of \eqref{RHSversion2}. Thus,
\begin{equation}\label{KGM''c}
  L^{\pm} \psi^{\pm} = \half A_{0} \bigl(\psi^{\pm} + e^{\pm 2 i tc^{2}} 
  \psi^{\mp}\bigr) \pm e^{\pm itc^{2}} \widetilde R,
\end{equation}
where
\begin{equation}\label{Rtilde}
  \widetilde R = \half M^{-1} \left\{ - 2ic \A \cdot \nabla \phi + 
  A_{0}(M-c^{2}) (\phi^{+} - \phi^{-}) + \A^{2} \phi \right\}.
\end{equation}

We are now ready to prove Theorem \ref{Thm1}.

\subsection{Spacetime norms}\label{Norms}

Define, for $0 \le T < \infty$,
\begin{align*}
  X_{T} &= \LpdotHslocal{\A(t)}{\infty}{1}
  + \tfrac{1}{c} \mixednormlocal{\partial_{t} \A(t)}{\infty}{2}
  + c \mixednormlocal{\square \A}{1}{2},
  \\
  Y_{T}^{\pm} &= \energylocal{\psi^{\pm}} + \LpHslocal{L^{\pm} 
  \psi^{\pm}}{1}{1},
  \\
  Z_{T}^{\pm} &= \mixednormlocal{\psi_{l}^{\pm}}{\infty}{2} + 
  \mixednormlocal{\psi_{l}^{\pm}}{2}{6},
\end{align*}
where $\psi_{l}^{\pm}$ is defined by \eqref{LowHigh}. Then set 
$Y_{T} = Y_{T}^{+} + Y_{T}^{-}$ and $Z_{T} = Z_{T}^{+} + Z_{T}^{-}$.
From the regularity properties \eqref{KMsolutionReg1}, 
\eqref{KMsolutionReg2} and \eqref{KMsolutionReg3} of 
$(A_{0},\A,\phi)$, it follows that
\begin{equation}\label{FiniteNorms}
  X_{T}, Y_{T}, Z_{T} < \infty.
\end{equation}
We prove this in section \ref{FiniteNormsProof} below.
Thus, $X_{T}, Y_{T}$ and $Z_{T}$ depend continuously on $T$. They also
depend on $c$, not only through the explicit appearance of $c$ in the 
definitions, but also through the implicit dependence of $\A$ 
and $\psi^{\pm}$ on $c$.

We claim that the assumptions on the data imply
\begin{equation}\label{InitialXYZbound}
  X_{0}, Y_{0}, Z_{0} = O(1)
\end{equation}
as $c \to \infty$. Obviously, \eqref{AdataBound} implies $X_{0} = O(1)$, 
and to bound $Y_{0}$ and $Z_{0}$ it suffices to check that
$\Sobnorm{\psi^{\pm}}{1} = O(1)$ at $t = 0$. But using 
\eqref{PhiSplitting} and Lemma \ref{Mestimates}\eqref{M1},
$$
  \Sobnorm{\psi^{\pm}}{1} = \Sobnorm{\phi^{\pm}}{1}
  \le \Sobnorm{\phi}{1} + \Sobnorm{M^{-1} \partial_{t} \phi}{1}
  + \tfrac{1}{c} \twonorm{A_{0} \phi}{},
$$
and by \eqref{PhidataBound} and \eqref{InitialEnergyBoundB}, the right hand 
side is $O(1)$ at $t = 0$.

Our main task will be to show that \eqref{InitialXYZbound} persists, 
i.e., for every $T < \infty$,
\begin{equation}\label{XYZbound}
  X_{T}, Y_{T}, Z_{T} = O(1)
\end{equation}
as $c \to \infty$.
In fact, we will prove \eqref{XYZbound} for a time $T = T_{0} > 0$
which only depends on the size of the global-in-time bound \eqref{L2bound}.
Then by iterating this argument we get \eqref{XYZbound} for every 
finite time $T$, since we can decompose $[0,T]$ into almost disjoint 
subintervals of length at most $T_{0}$.

Once \eqref{XYZbound} has been proved, the local-in-time bounds in Theorem 
\ref{Thm1} follow easily, as we demonstrate in section 
\ref{ConclusionThm1}.

\subsection{Main estimates and bootstrap argument}

Here we prove \eqref{XYZbound} for a time $T = T_{0} > 0$
which only depends on the size of \eqref{L2bound}.
Using a bootstrap argument, we reduce this to proving (recall the notational 
conventions made in section \ref{MainResults})
\begin{align}
  \label{Xineq}
  X_{T} &\lesssim X_{0} + \sigma(T) Y_{T}^{2}
  + \tfrac{\sigma(T)}{c} P(X_{T}) P(Y_{T}),
  \\
  \label{Yineq}
  Y_{T} &\lesssim Y_{0} + \sigma(T) Z_{T}^{2} Y_{T} + 
  \tfrac{\sigma(T)}{c^{1/2}} P(X_{T}) P(Y_{T}),
  \\
  \label{Zineq}
  Z_{T} &\lesssim Z_{0} + \sigma(T) \Bigl(\sum\nolimits_{\pm}
  \mixednorm{\psi^{\pm}}{\infty}{2}^{2} \Bigr) Z_{T}
  + \tfrac{\sigma(T)}{c^{1/2}} P(X_{T}) P(Y_{T}),
\end{align}
for, say, $0 \le T \le 1$ and $c \ge 1$.

Indeed, assuming these inequalities hold, first observe that 
\eqref{Zineq} implies
\begin{equation}\label{Zineq'}
  Z_{T} \lesssim Z_{0} + \tfrac{\sigma(T)}{c^{1/2}} P(X_{T}) P(Y_{T})
  \quad \text{for} \quad 0 \le T \le T_{0},
\end{equation}
for some $T_{0} > 0$ which only depends on \eqref{L2bound}. Plugging this 
into \eqref{Yineq} gives
$$
  Y_{T} \lesssim Y_{0} + \sigma(T) Z_{0}^{2} Y_{T} + 
  \tfrac{\sigma(T)}{c^{1/2}} P(X_{T}) P(Y_{T})
  \quad \text{for} \quad 0 \le T \le T_{0}.
$$
Thus, making $T_{0}$ smaller if necessary, but still depending only 
on \eqref{L2bound}, we get
\begin{equation}\label{Yineq'}
  Y_{T} \lesssim Y_{0} + 
  \tfrac{\sigma(T)}{c^{1/2}} P(X_{T}) P(Y_{T})
  \quad \text{for} \quad 0 \le T \le T_{0}.
\end{equation}
Inserting this into the second term on the right hand side of \eqref{Xineq} gives
\begin{equation}\label{Xineq'}
  X_{T} \lesssim X_{0} + Y_{0}^{2}
  + \tfrac{\sigma(T)}{c} P(Y_{T}) P(X_{T})
  \quad \text{for} \quad 0 \le T \le T_{0}.
\end{equation}
Adding up \eqref{Xineq'} and \eqref{Yineq'} gives
\begin{equation}\label{BootStrapIneq}
  f(T) \le P(f(0)) + \tfrac{\sigma(T)}{c^{1/2}} Q(f(T)) f(T)
  \quad \text{for} \quad 0 \le T \le T_{0},
\end{equation}
where $f(T) = X_{T} + Y_{T}$ depends continuously on $T$ and $Q$ is a polynomial.
We claim that \eqref{BootStrapIneq} implies
\begin{equation}\label{fClaim}
  f(T) \le 2 P(f(0)) \quad \text{for} \quad 0 \le T \le T_{0}
\end{equation}
if $c$ is sufficiently large (depending on $f(0)$). In view of 
\eqref{Zineq'} and \eqref{InitialXYZbound}, this implies 
\eqref{XYZbound} for $T \le T_{0}$.

Let us prove \eqref{fClaim}. We consider two cases: $f(0) = 0$ or $f(0) > 0$.
If $f(0) = 0$, then letting $c \to \infty$ in \eqref{BootStrapIneq} would lead
to a contradiction unless $f(T) = 0$ for $0 \le T \le T_0$.
On the other hand, if $f(0) > 0$, then we claim that \eqref{fClaim} must hold with
strict inequality if $c$ is sufficiently large. For if not, then by continuity,
$f(T) = 2 P(f(0))$ for some $0 \le T \le T_{0}$, which by \eqref{BootStrapIneq} implies
$$
  P(f(0)) \le \tfrac{\sigma(T)}{c^{1/2}} Q[2P(f(0))] 2 P(f(0)).
$$
Dividing by $P(f(0))$ gives
$$
  \sigma(T) \ge \frac{c^{1/2}}{2Q[2P(f(0))]},
$$
but this fails for sufficiently large $c$.

Thus, we have reduced \eqref{XYZbound} to proving \eqref{Xineq}--\eqref{Zineq}.
To do this, we will use energy estimates and the spacetime estimates proved in
section \ref{SpacetimeEstimates}. Let us turn to the details.
We start by proving some estimates for the elliptic variable 
$A_{0}$.

\subsubsection{Estimates for $A_{0}$}\label{A0Estimates}

Our aim here is to prove:

\begin{lemma}\label{A0Prop}
Let $1 \le q < 2$. Then
\begin{enumerate}
  \item\label{A0Prop1}
  $\mixednormlocal{\Delta A_{0}}{q}{r}
  \lesssim
  \sigma(T) Z_{T}^{2} + \tfrac{\sigma(T)}{c^{1/2}} Y_{T}^{2}$
  if $1 \le r \le \tfrac{3}{2}^{+}$.
  \item\label{A0Prop2}
  $\mixednormlocal{\Delta A_{0}}{q}{r}
  \lesssim
  \sigma(T) \left( \sum\nolimits_{\pm}
  \mixednormlocal{\psi^{\pm}}{\infty}{2} \right) Z_{T}
  + \tfrac{\sigma(T)}{c^{1/2}} Y_{T}^{2}$
  if $1 \le r \le \tfrac{3}{2}$.
  \item\label{A0Prop3}
  $\mixednormlocal{\Delta A_{0}}{\infty}{r}
  \lesssim Y_{T}^{2}$
  if $1 \le r \le \tfrac{3}{2}$.
\end{enumerate}
\end{lemma}

Expanding RHS\eqref{KGM'a} in terms of $\phi^{\pm} = e^{\mp 
itc^{2}} \psi^{\pm}$, we reduce this to proving the same estimates for
$\norm{I}_{L_{t}^{q}([0,T])}$ where
$I(t) = \bigLxpnorm{\psi^{\pm} \tfrac{M}{c^{2}} \overline{\psi^{\pm}}}{r}$
and the $\pm$ signs are independent.
Expanding $\psi^{\pm} = \psi^{\pm}_{l} + \psi^{\pm}_{h}$
as in \eqref{LowHigh} gives
$$
  I
  \le I_{l,l} + I_{l,h} + I_{h,l} + I_{h,h}
  \quad \text{where} \quad
  I_{\cdot,\cdot} = \Lxpnorm{\psi^{\pm}_{(\cdot)}\tfrac{M}{c^{2}}
  \overline{\psi^{\pm}_{(\cdot)}}}{r}.
$$

\paragraph{The case $r = \tfrac{3}{2}^{+}$.}

By H\"older's inequality, Lemma \ref{LowHighLemma}\eqref{HL1} and 
$L^{p}$ interpolation,
\begin{equation}\label{IlowlowFixedTime}
  I_{l,l}
  \le \Lxpnorm{\psi^{\pm}_{l}}{2^{+}} \Lxpnorm{\psi^{\pm}_{l}}{6}
  \lesssim \Lxpnorm{\psi^{\pm}_{l}}{2}^{1-\delta} 
  \Lxpnorm{\psi^{\pm}_{l}}{6}^{1+\delta}
\end{equation}
for some $\delta > 0$. Since $q < 2$, and $\delta \to 0$ as $r \to 3/2$,
we will have $\tfrac{1}{q} - \tfrac{1+\delta}{2} > 0$ if $r$ is close 
enough to $3/2$. Applying H\"older's inequality in $t$ then yields
$$
  \norm{I_{l,l}}_{L_{t}^{q}([0,T])}
  \lesssim \mixednormlocal{\psi^{\pm}_{l}}{\infty}{2}^{1-\delta}
  T^{\frac{1}{q} - \frac{1+\delta}{2}} 
  \mixednormlocal{\psi^{\pm}_{l}}{2}{6}^{1+\delta}
  \lesssim \sigma(T) Z_{T}^{2}
$$
as desired. Next, by H\"older's inequality 
and Lemma \ref{LowHighLemma},
\begin{align}
  \label{Ilowhigh}
  I_{l,h} &\lesssim
  \Lxpnorm{\psi^{\pm}_{l}}{6^{+}}
  \Lxpnorm{\tfrac{M}{c^{2}}\psi^{\pm}_{h}}{2}
  \lesssim
  \Lxpnorm{\psi^{\pm}_{l}}{6^{+}}
  \tfrac{1}{c} \Sobnorm{\psi^{\pm}}{1},
  \\
  \label{Ihighlow}
  I_{h,l} &\lesssim
  \Lxpnorm{\psi^{\pm}_{h}}{2}
  \Lxpnorm{\tfrac{M}{c^{2}}\psi^{\pm}_{l}}{6^{+}}
  \lesssim
  \tfrac{1}{c} \Sobnorm{\psi^{\pm}}{1}
  \Lxpnorm{\psi^{\pm}_{l}}{6^{+}},
  \\
  \label{Ihighhigh}
  I_{h,h} &\lesssim 
  \Lxpnorm{\psi^{\pm}_{h}}{6^{+}}
  \Lxpnorm{\tfrac{M}{c^{2}}\psi^{\pm}_{h}}{2}
  \lesssim
  \Lxpnorm{\psi^{\pm}_{h}}{6^{+}}
  \tfrac{1}{c} \Sobnorm{\psi^{\pm}}{1}.
\end{align}
Since $\Lxpnorm{\psi^{\pm}_{l}}{6^{+}} \lesssim c^{\varepsilon} 
\Sobnorm{\psi^{\pm}}{1}$ by Sobolev embedding and Lemma 
\ref{LowHighLemma}\eqref{HL2}, it follows that
$I_{l,h}, I_{h,l} \lesssim c^{\varepsilon-1} Y_{T}^{2}$,
whence
$$
  \norm{I_{l,h}}_{L_{t}^{q}([0,T])}, \norm{I_{l,h}}_{L_{t}^{q}([0,T])}
  \lesssim \tfrac{T^{1/q}}{c^{1/2}} Y_{T}^{2}
$$
as desired. Finally, to control $I_{h,h}$, we have to use Strichartz estimates.
Applying H\"older's inequality in $t$ to \eqref{Ihighhigh} gives
$$
  \norm{I_{h,h}}_{L_{t}^{q}([0,T])}
  \lesssim T^{\tfrac{1}{q} - \tfrac{1}{p}}
  \mixednormlocal{\psi^{\pm}}{p}{6^{+}}
  \tfrac{1}{c} \energylocal{\psi^{\pm}}.
$$
Choosing $p = 3^{-}$ so that $(p,6^{+})$ is sharp wave 
admissible, we have
$$
  \mixednormlocal{\psi^{\pm}}{p}{6^{+}} \lesssim Y_{T}
$$
by the corollary to Proposition \ref{WaveProp}. This proves part 
\eqref{A0Prop1} of Lemma \ref{A0Prop}.

\paragraph{The case $r \le \tfrac{3}{2}$.}
This is similar, but simpler. Instead of \eqref{IlowlowFixedTime}, we have
$$
  I_{l,l}
  \le \Lxpnorm{\psi^{\pm}_{l}}{2} \Lxpnorm{\psi^{\pm}_{l}}{p}
  \le \Lxpnorm{\psi^{\pm}_{l}}{2}^{1+\delta} 
  \Lxpnorm{\psi^{\pm}_{l}}{6}^{1-\delta}
$$
for some $2 \le p \le 6$ and $0 \le \delta \le 1$. Thus,
\begin{align*}
  \norm{I_{l,l}}_{L_{t}^{\infty}([0,T])}
  &\lesssim 
  \energylocal{\psi^{\pm}}^{2} \le Y_{T}^{2},
  \\
  \norm{I_{l,l}}_{L_{t}^{q}([0,T])}
  &\lesssim
  \mixednormlocal{\psi^{\pm}_{l}}{\infty}{2}^{1+\delta}
  T^{\frac{1}{q} - \frac{1-\delta}{2}} 
  \mixednormlocal{\psi^{\pm}_{l}}{2}{6}^{1-\delta},
\end{align*}
as desired for parts \eqref{A0Prop3} and \eqref{A0Prop2}, respectively, 
of Lemma \ref{A0Prop}. Next, since the estimates 
\eqref{Ilowhigh}--\eqref{Ihighhigh} now hold with $6^{+}$ replaced by 
some $2 \le p \le 6$, we have
\begin{equation}\label{errorBounds}
  I_{l,h}, I_{h,l}, I_{h,h} \lesssim \tfrac{1}{c} 
  \energylocal{\psi^{\pm}}^{2} \le \tfrac{1}{c} Y_{T}^{2}.
\end{equation}
by Sobolev embedding. This concludes the proof of Lemma \ref{A0Prop}.

\subsubsection{Estimate for $X_{T}$}

Here we prove \eqref{Xineq}. By the energy inequality for \eqref{IVPforWave}, which reads
$$
  \LpdotHslocal{u}{\infty}{1}
  + \tfrac{1}{c} \mixednormlocal{\partial_{t} u}{\infty}{2}
  \le \Sobdotnorm{f_{0}}{1} + \tfrac{1}{c} \twonorm{f_{1}}{}
  + c \mixednormlocal{F}{1}{2},
$$
we have $X_{T} \le X_{0} + c \mixednormlocal{\square \A}{1}{2}$,
so in view of \eqref{KGM'b} it suffices to prove
\begin{gather}
  \label{XproofA}
  \mixednormlocal{\Proj\bigl( \phi \overline{\nabla \phi} \bigr)}{1}{2}
  \lesssim \sqrt{T} Y_{T}^{2},
  \\
  \label{XproofB}
  \bigmixednormlocal{\abs{\phi}^{2} \A}{1}{2}
  \lesssim
  T Y_T^2 X_T.
\end{gather}
To do this, substitute
\begin{equation}\label{PhiExpansion}
  \phi = \phi^{+} + \phi^{-} = e^{-itc^{2}} \psi^{+} + e^{+itc^{2}} 
  \psi^{-}
\end{equation}
in the left hand sides, and expand.

Applying H\"older's inequality in $t$, we conclude that
LHS\eqref{XproofA} is bounded by a sum of terms
$\sqrt{T} \bigtwonorm{\Proj\bigl( \psi^{\pm} \overline{\nabla \psi^{\pm}} 
\bigr)}{(S_{T})}$, where the signs are independent.
Noting the identity $\overline{L^{\pm} u} = - L^{\mp} \overline u$, we 
apply the corollary to Proposition \ref{EasyNullFormProp} to conclude
that \eqref{XproofA} holds.

Next, LHS\eqref{XproofB} is bounded by a sum of terms $T 
\bigmixednormlocal{\psi^{\pm} \overline{ \psi^{\pm}} \A}{\infty}{2}$,
with independent signs. Applying the inequality \eqref{L2trilinear} gives \eqref{XproofB}.

\subsubsection{Estimate for $Y_{T}$}

Here we prove \eqref{Yineq}. In fact, we only prove $Y_{T}^{+} \lesssim$ 
RHS\eqref{Yineq}, as the proof for $Y_{T}^{-}$ is the same. For the solution of 
\eqref{IVPforL} we have, in view of the formula \eqref{Lsolution},
$$
  \LpHslocal{u}{\infty}{1} \le \Sobnorm{f}{1} + \LpHslocal{F}{1}{1}.
$$  
Thus, $Y_{T}^{+} \lesssim Y_{0}^{+} + \LpHslocal{L^{+}\psi^{+}}{1}{1}$,
so in view of \eqref{KGM''c} and Lemma \ref{Mestimates}\eqref{M1},
it suffices to prove
\begin{align}
  \label{YproofA}
  \LpHslocal{A_{0} \psi^{\pm}}{1}{1}
  &\lesssim \sigma(T) Z_{T}^{2} Y_{T} + \tfrac{\sigma(T)}{c^{1/2}} Y_{T}^{3},
  \\
  \label{YproofB}
  \mixednormlocal{\A \cdot \nabla \phi}{1}{2}
  &\lesssim \sqrt{\tfrac{T}{c}} X_{T} Y_{T},
  \\
  \label{YproofC}
  \mixednormlocal{\A^{2} \phi}{1}{2}
  &\lesssim T X_{T}^{2} Y_{T},
  \\
  \label{YproofD}
  \mixednormlocal{A_{0} \tfrac{M-c^{2}}{c} \psi^{\pm}}{1}{2}
  &\lesssim \sigma(T) Z_{T}^{2} Y_{T} + \tfrac{\sigma(T)}{c^{1/2}} Y_{T}^{3}.
\end{align}

To prove \eqref{YproofB} and \eqref{YproofC}, expand using 
\eqref{PhiExpansion}, and apply, respectively, the corollary to 
Proposition \ref{HardNullFormProp} and inequality \eqref{L2trilinear}.

Next, observe that by the product rule for derivatives, H\"older's 
inequality and the Sobolev embedding \eqref{L6Sobolev},
the LHS\eqref{YproofA} is dominated by
\begin{equation}\label{YproofE}
  \int_{0}^{T} \left( 
  \Lpnorm{\nabla A_{0}}{3} + \inftynorm{A_{0}}{} \right) \, dt 
  \, \energylocal{\psi^{\pm}},
\end{equation}
and by H\"older's inequality and Lemma \ref{Mestimates}\eqref{M3}, 
LHS\eqref{YproofD} is also $\lesssim$ \eqref{YproofE}.
Thus, it is enough to show
\begin{equation}\label{A0Est}
  \int_{0}^{T} \left( 
  \Lpnorm{\nabla A_{0}}{3} + \inftynorm{A_{0}}{} \right) \, dt
  \lesssim
  \sigma(T) Z_{T}^{2} + \tfrac{\sigma(T)}{c^{1/2}} Y_{T}^{2}.
\end{equation}
But in view of Lemma \ref{EllipticLemma}, this follows from
Lemma \ref{A0Prop} (section \ref{A0Estimates}).

\subsubsection{Estimate for $Z_{T}$}\label{Zestimate}

Here we prove \eqref{Zineq}. We only prove $Z_{T}^{+} \lesssim$ 
RHS\eqref{Zineq}, as $Z_{T}^{-}$ can be treated in the same way.
Our argument is reminiscent of that used in section \ref{PSbounds}
to prove the $L^2$ well-posedness of the Schr\"odinger-Poisson system.
Thus, we apply the Strichartz estimate in Proposition \ref{SchrodingerProp} 
with $(q,r) = (2,6)$ and $(\widetilde q', \widetilde r') = (1,2)$ or
$(1^{+},2^{-})$. More precisely, using \eqref{KGM''c} we write
$L^{+} \psi^{+} = F + G$, where
\begin{align*}
  2 F &=  A_{0} \bigl(\psi^{+} + e^{2 i tc^{2}} \psi^{-}\bigr),
  \\
  2 G &= e^{itc^{2}} M^{-1}\left\{ - 2ic \A \cdot \nabla \phi + 
  A_{0}(M-c^{2}) (\phi^{+} - \phi^{-}) + \A^{2} \phi \right\}.
\end{align*}
Then by Proposition \ref{SchrodingerProp},
$$
  Z_{T}^{+} \lesssim Z_{0}^{+} + \mixednormlocal{F}{a}{b}
  + \mixednormlocal{G}{1}{2},
$$
where $\tfrac{1}{a} + \tfrac{\varepsilon}{2} = 1$ and
$\tfrac{1}{b} + \tfrac{1}{2} - \tfrac{\varepsilon}{3} = 1$ for some
 sufficiently small $\varepsilon > 0$.
Thus, it suffices to prove
\begin{gather}
  \label{Fest}
  \mixednormlocal{F}{a}{b} \lesssim
  \sigma(T) \left( \sum\nolimits_{\pm} 
  \mixednormlocal{\psi^{\pm}}{\infty}{2}^{2} \right) Z_{T}
  + \tfrac{\sigma(T)}{c^{1/2}} Y_{T}^{3},
  \\
  \label{Gest}
  \mixednormlocal{G}{1}{2} \lesssim 
  \tfrac{\sigma(T)}{c} P(X_{T}) P(Y_{T}).
\end{gather}
First, write
$$
  \mixednormlocal{F}{a}{b}
  \le \mixednormlocal{A_{0}}{a}{3/\varepsilon} \sum\nolimits_{\pm} 
  \mixednormlocal{\psi^{\pm}}{\infty}{2}.
$$
Since 
$\mixednormlocal{A_{0}}{a}{3/\varepsilon} \lesssim 
\mixednormlocal{\Delta A_{0}}{a}{(3/2)^{-}}$ by Sobolev embedding,
\eqref{Fest} then follows from Lemma \ref{A0Prop}\eqref{A0Prop2}.
Next, observe that \eqref{Gest} follows from Lemma \ref{Mestimates}\eqref{M2}
and the estimates \eqref{YproofB}, \eqref{YproofC} and 
$$
  \mixednormlocal{A_{0} \tfrac{M-c^{2}}{c} \psi^{\pm}}{1}{2}
  \lesssim \sigma(T) Y_{T}^{3}.
$$
This last inequality follows from \eqref{YproofD} and the fact that
\begin{equation}\label{ZcontrolledbyY}
  Z_{T} \lesssim (1 + T^{1/2}) \sum\nolimits_{\pm} 
  \energylocal{\psi^{\pm}}
  \le (1 + T^{1/2}) Y_{T},
\end{equation}
where we use \eqref{L6Sobolev} to get the first inequality.

\subsubsection{Finiteness of norms}\label{FiniteNormsProof}

Here we prove the claim made earlier, that the regularity properties 
of $(A_{0},\A,\phi)$ imply \eqref{FiniteNorms}.
First, $X_{T} < \infty$ follows directly from \eqref{KMsolutionReg1}
and \eqref{KMsolutionReg2}. 
Next, using the definition \eqref{PhiSplitting}, Lemma \ref{Mestimates}\eqref{M1} and
\eqref{L2bilinear}, we conclude from \eqref{KMsolutionReg1} that
\begin{equation}\label{PhipmFinite}
  \phi^{\pm} \in C(\R;H^{1}).
\end{equation}
In view of \eqref{ZcontrolledbyY}, this implies $Z_{T} < \infty$.
Moreover, it reduces $Y_{T} < \infty$ to showing that
$$
  \LpHslocal{L^{\pm} \psi^{\pm}}{1}{1} < \infty.
$$
But the latter reduces to proving that the left hand sides of 
\eqref{YproofA}--\eqref{YproofD} are finite.
First recall that LHS\eqref{YproofA} and LHS\eqref{YproofD} are bounded by \eqref{YproofE}, 
which is finite by \eqref{KMsolutionReg3} and \eqref{PhipmFinite}.
Next, using \cite[Proposition 2.2]{KM} instead of the corollary to Proposition
\ref{HardNullFormProp}, one finds that LHS\eqref{YproofB} is controlled
by \eqref{KMsolutionReg2} and the norms of the initial data \eqref{KGMdata}. 
Finally,  LHS\eqref{YproofC} $<\infty$ by \eqref{KMsolutionReg1}, if we
use \eqref{L2trilinear}.

\subsection{Conclusion of proof of Theorem 
\ref{Thm1}}\label{ConclusionThm1}

We conclude by showing that \eqref{XYZbound} implies the 
local-in-time bounds in Theorem \ref{Thm1}.
By the definitions of $X_{T}, Y_{T}$ and $Z_{T}$, 
it is obvious that they control the norms in parts 
\eqref{Thm1A}--\eqref{Thm1C} in Theorem \ref{Thm1}.
The bound \eqref{Thm1F} reduces to
Lemma \ref{A0Prop}\eqref{A0Prop3} via Sobolev embedding and 
\eqref{RieszCorollary}.
To prove part \eqref{Thm1H}, use Lemma 
\ref{Mestimates}\eqref{M1} and \eqref{L2bilinear} to get
$$
  \Sobnorm{M^{-1}(A_{0}\phi)}{1} \lesssim \tfrac{1}{c}
  \twonorm{\nabla A_{0}}{} \Sobnorm{\phi}{1}^{2}
$$
for each $t$. Then use the bounds in parts \eqref{Thm1A} and 
\eqref{Thm1F}. It only remains to prove the bound for $R$ in
part \eqref{Thm1G} of Theorem \ref{Thm1}. By Lemma 
\ref{Mestimates}\eqref{M1}, this reduces to $\tfrac{1}{c} 
\mixednormlocal{M R}{1}{2} = O(c^{-1/2})$.
Recalling the definition \eqref{Remainder} of $R$
and the estimates \eqref{YproofB} and \eqref{YproofC}, we see that it 
suffices to prove
$$
  \mixednormlocal{[ A_{0}, M-c^{2} ] \psi^{\pm}}{1}{2} = O(c^{1/2}).
$$
To do this, expand the RHS\eqref{KGM'a} using the frequency 
decomposition \eqref{LowHigh}, as in the proof of Lemma 
\ref{A0Prop}, and write $A_{0} = A_{0}' + A_{0}''$, where $A_{0}'$ 
corresponds to terms of the type
$\psi^{\pm}_{l} \tfrac{M}{c^{2}} \overline{\psi^{\pm}_{l}}$, i.e., 
both factors are at low frequency, and $A_{0}''$ corresponds to terms 
where at least one factor has high frequency. Let us consider first
$$
  \mixednormlocal{[ A_{0}'', M-c^{2} ] \psi^{\pm}}{1}{2}.
$$
Here we do not need the commutator structure, so we simply use 
\eqref{Commutator} and Lemma \ref{Mestimates}\eqref{M3} to dominate it 
by
$$
  c \left( \mixednormlocal{\nabla A_{0}''}{1}{3} + 
  \mixednormlocal{A_{0}''}{1}{\infty} \right) \energylocal{\psi^{\pm}}
$$
In view of Lemma \ref{EllipticLemma}, it therefore suffices to check
$$
  \mixednormlocal{\Delta A_{0}''}{1}{(3/2)^{\pm}} = O(c^{-1/2}),
$$
but this is clear from the proof of Lemma \ref{A0Prop}, since for
$A_{0}''$ there is no term $I_{l,l}$. It remains to prove
$$
  \mixednormlocal{[ A_{0}', M-c^{2} ] \psi^{\pm}}{1}{2} = O(c^{1/2}).
$$
In fact, applying the following lemma with $f = \psi^{\pm}_{l}, g = 
\tfrac{M}{c^{2}} \psi^{\pm}_{l}$ and $h = \psi^{\pm}$, and using 
Lemma \ref{LowHighLemma}, parts \eqref{HL1} and \eqref{HL2}, gives
$\mixednormlocal{[ A_{0}', M-c^{2} ] \psi^{\pm}}{1}{2}
\lesssim c^{\varepsilon} Y_{T}^{3}$ for $\varepsilon > 0$ arbitrarily small.

\begin{lemma}\label{CommutatorLemma}
Define
$$
  T(f,g,h) = (M-c^{2}) \Bigl[ (-\Delta)^{-1}(fg) \cdot h \Bigr]
  - (-\Delta)^{-1}(fg) \cdot (M-c^{2}) h.
$$
Then the estimate
$\twonorm{T(f,g,h)}{} \lesssim \Sobnorm{f}{1^{+}} \Sobnorm{g}{1^{+}} 
\Sobnorm{h}{1}$ holds on $\R^{3}$.
\end{lemma}

\begin{proof} The Fourier symbol of $T$ is
$$
  \sigma(\xi,\eta,\zeta) = \frac{1}{\abs{\xi+\eta}^{2}} \bigl[ 
  h_{c}(\xi+\eta+\zeta) - h_{c}(\zeta) \bigr],
$$
where $h_{c}$ is given by \eqref{hc}. We claim that
\begin{equation}\label{SymbolClaim}
  \abs{\sigma(\xi,\eta,\zeta)} \lesssim 1 + 
  \frac{\abs{\zeta}}{\abs{\xi+\eta}}.
\end{equation}
Since we may assume that $\widehat f, \widehat g, \widehat h \ge 0$, 
this would imply
$$
  \twonorm{T(f,g,h)}{} \lesssim \twonorm{fgh}{} + 
  \twonorm{\D^{-1}(fg) \cdot \D h}{}.
$$
The first term on the right hand side is covered by 
\eqref{L2trilinear}, the second term is $\le$
$$
  \inftynorm{\D^{-1}(fg)}{} \Sobnorm{h}{1} \lesssim 
  \Sobnorm{fg}{(1/2)^{+}} \Sobnorm{h}{1}.
$$
Since $\Sobnorm{fg}{(1/2)^{+}} \lesssim \Sobnorm{f}{1^{+}} 
\Sobnorm{g}{1^{+}}$, we get the desired estimate.

Thus, it only remains to prove \eqref{SymbolClaim}, which clearly 
reduces to
$$
  \abs{h_{c}(\xi+\eta) - h_{c}(\eta)} \lesssim \abs{\xi} \bigl( 
  \abs{\xi} + \abs{\eta} \bigr).
$$
By the mean value theorem, this reduces to checking
\begin{equation}\label{Gradientbound}
  \abs{\nabla h_{c}(\xi)} \lesssim \abs{\xi}.
\end{equation}
But writing $h = h_{1}$, we have $h_{c}(\xi) = c^{2} h(\xi/c)$, so
$$
  \nabla h_{c}(\xi) = c (\nabla h)(\xi/c),
$$
and we know from Lemma \ref{alphaLemma} that $\abs{\nabla h(\xi)} 
\lesssim \abs{\xi}$ for all $\xi$.
\end{proof}

\section{Proof of $H^{1}$ convergence}\label{H1Convergence}

Here we prove Theorem \ref{Thm2}. Thus, we assume that the hypotheses of 
the theorem are satisfied, with one modification: As noted in section 
\ref{MainResults}, in view of the bound in Theorem 
\ref{Thm1}\eqref{Thm1H} we may use the definition 
\eqref{PhiSplitting} instead of \eqref{FreeSplitting}. The 
equations \eqref{KGM'} are therefore satisfied.

We first prove \eqref{Convergence}.
Clearly, it is enough to show that given $0 < T < \infty$, there exist
constants $B, \varepsilon > 0$ such that for every time interval
$I = [t_{0},t_{1}] \subset [0,T]$,
\begin{equation}\label{ReducedConvergence}
  f(I) \lesssim
  \sum\nolimits_{\pm} \Sobnorm{\psi^{\pm}(t_{0}) - v^{\pm}(t_{0})}{1}
  + B \abs{I}^{\varepsilon} f(I) + o(1),
\end{equation}
where $f(I) = \sum\nolimits_{\pm} \norm{\psi^{\pm} -
v^{\pm}}_{L_{t}^{\infty}H^{1}(I \times \R^{3})}$. In fact, $B$ and 
$\varepsilon$ depend only on the bounds in Theorem \ref{Thm1} and
Lemma \ref{PSlemma}.

Without loss of generality, we assume $I = [0,T]$, and we choose the 
plus sign on the left hand side of \eqref{ReducedConvergence}.
By the formula \eqref{Lsolution} for the initial value problem 
\eqref{IVPforL}, and the corresponding formula for the Schr\"odinger 
equation,\begin{align*}
  \psi^{+}(t) &= U^{+}(t) \psi^{+}(0) + \int_{0}^{t} 
  U^{+}(t-s)(L^{+}\psi^{+})(s) \, ds,
  \\
  v^{+}(t) &= V^{+}(t) v^{+}(0) + \int_{0}^{t} 
  V^{+}(t-s)(uv^{+})(s) \, ds,
\end{align*}
where $U^{+}$ and $V^{+}$ are given by \eqref{Lpropagator} and 
\eqref{SchrodingerPropagator}. Thus,
\begin{align*}
  \psi^{+}(t) - v^{+}(t) &=
  U^{+}(t) [ \psi^{+}(0) - v^{+}(0) ] + [ U^{+}(t) - V^{+}(t) ] 
  v^{+}(0)
  \\
  &\quad + \int_{0}^{t} U^{+}(t-s) [L^{+}\psi^{+} - uv^{+}](s) \, ds
  \\
  &\quad + \int_{0}^{t} [ U^{+}(t-s) - V^{+}(t-s) ] (uv^{+})(s) \, ds
  \\
  &= I_{1} + \cdots + I_{4}.
\end{align*}
Now, $(U^{+}(t) f)\,\widehat\,(\xi) = e^{-ith_{c}(\xi)} \widehat f(\xi)$, 
with $h_{c}$ given by \eqref{hc}, and $(V^{+}(t) f)\,\widehat\,(\xi) = 
e^{it\abs{\xi}^{2}/2} \widehat f(\xi)$. Using Plancherel's theorem, 
it is therefore clear that
$$
  \energylocal{I_{1}} \le \Sobnorm{\psi^{+}(0) - v^{+}(0)}{1}.
$$
Moreover,
\begin{align*}
  \energylocal{I_{2}}
  &\le \twonorm{ \sup_{0 \le t \le T} \abs{ e^{-ith_{c}(\xi)} - 
  e^{it\abs{\xi}^{2}/2} } (1+\abs{\xi}^{2})^{1/2} \bigabs{\widehat{v^{+}}(0,\xi)} 
  }{_{\xi}}
  \\
  &\lesssim
  \twonorm{ \min\left\{ 1, T \abs{h_{c}(\xi) - \tfrac{\abs{\xi}^{2}}{2}}
  \right\}
  (1+\abs{\xi}^{2})^{1/2} \bigabs{\widehat{v^{+}}(0,\xi)} }{_{\xi}},
\end{align*}
and the latter $\to 0$ as $c \to \infty$, by the dominated convergence theorem.
Similarly, with $F = uv^{+}$,
$$
  \energylocal{I_{4}} \lesssim
  \int_{0}^{T} \twonorm{ \min\left\{ 1, T \abs{h_{c}(\xi) - \tfrac{\abs{\xi}^{2}}{2}}
  \right\}
  (1+\abs{\xi}^{2})^{1/2} \bigabs{\widehat F(s,\xi)} }{_{\xi}} \, ds,
$$
and this $\to 0$ as $c \to \infty$ by dominated convergence, because
$\energylocal{F} < \infty$. To prove the latter, note that (cf. 
\eqref{YproofE})
$$
  \Sobnorm{uv^{+}}{1} \lesssim
  \left(\Lxpnorm{\nabla u}{3} + \Lxpnorm{u}{\infty} \right)
  \Sobnorm{v^{+}}{1}
$$
and
\begin{equation}\label{uBound}
  \Lxpnorm{\nabla u}{3} + \Lxpnorm{u}{\infty}
  \lesssim
  \Lxpnorm{\Delta u}{(3/2)^{-}} + \Lxpnorm{\Delta u}{(3/2)^{+}}
  \lesssim \sum\nolimits_{\pm} \Sobnorm{v^{\pm}}{1}^{2},
\end{equation}
where we used Lemma \ref{EllipticLemma} to get the first inequality, 
then H\"older's inequality and Sobolev embedding to get the second one.

It only remains to estimate $I_{3}$. Write
$$
  L^{+}\psi^{+} - uv^{+} = (A_{0}-u) \psi^{+} + u (\psi^{+}-v^{+})
  + e^{itc^{2}} R,
$$
where $R$ is given by \eqref{Remainder}. Correspondingly, we split 
$I_{3} = I_{3}' + I_{3}'' + I_{3}'''$.
First observe that
$$
  \energylocal{I_{3}'''} \le \LpHslocal{R}{1}{1} = O(c^{-1/2})
$$
by part \eqref{Thm1G} of Theorem \ref{Thm1}. Next, write
\begin{align*}
  \energylocal{I_{3}''} &\le \int_{0}^{T} \Sobnorm{u (\psi^{+}-v^{+})}{1} \, dt
  \\
  &\lesssim
  T
  \left(\mixednormlocal{\nabla u}{\infty}{3} + \inftynorm{u}{(S_{T})} \right)
  \energylocal{\psi^{+}-v^{+}}
\end{align*}
and recall \eqref{uBound}. Similarly,
\begin{align*}
  \energylocal{I_{3}'} \lesssim
  \left(\mixednormlocal{\nabla (A_{0}- u)}{1}{3} + 
  \mixednormlocal{A_{0} - u}{1}{\infty} \right)
  \energylocal{\psi^{+}},
\end{align*}
so in view of Lemma \ref{EllipticLemma}, to finish the proof of 
\eqref{ReducedConvergence} it only remains to show
\begin{equation}\label{Difference}
  \mixednormlocal{\Delta(A_{0} - u)}{1}{r}
  \le B \sigma(T) \sum\nolimits_{\pm} \energylocal{\psi^{\pm}-v^{\pm}}
  + o(1)
\end{equation}
for $1 \le r \le \tfrac{3}{2}^{+}$. To this end, observe that by
\eqref{DeltaA0Expansion}, \eqref{A0Error} and \eqref{PSa},
$$
  \Delta(A_{0} - u) = - \bigabs{\psi^{+}}^{2} + \bigabs{v^{+}}^{2}
  + \bigabs{\psi^{-}}^{2} - \bigabs{v^{-}}^{2} + \text{Error},
$$
where
\begin{equation}\label{Error}
  \Lxpnorm{\text{Error}}{r} \lesssim \sum \Lxpnorm{\psi^{\pm} 
  \left( \tfrac{M}{c^{2}} - 1 \right) \overline{\psi^{\pm}}}{r},
\end{equation}
and the sum is over all combinations of signs.

Since $\abs{\psi^{+}}^{2} - \abs{v^{+}}^{2} = \left( \abs{\psi^{+}} - 
\abs{v^{+}} \right) \left( \abs{\psi^{+}} + \abs{v^{+}} \right)$,
H\"older's inequality and Sobolev embedding yield
$$
  \Lxpnorm{\abs{\psi^{+}}^{2} - \abs{v^{+}}^{2}}{r}
  \lesssim \Sobnorm{\psi^{+} - v^{+}}{1} \left( 
  \Sobnorm{\psi^{+}}{1} + \Sobnorm{v^{+}}{1} \right).
$$
Integrating in time and using the bounds in Theorem \ref{Thm1} and 
Lemma \ref{PSlemma} then gives the first term on the right hand side 
of \eqref{Difference}.

It only remains to prove $\mixednormlocal{\text{Error}}{1}{r} 
\lesssim$ RHS\eqref{Difference}.
Splitting the right hand side of \eqref{Error}
as in the proof of Lemma \ref{A0Prop}, and using the 
estimates obtained there, we get
$$
  \mixednormlocal{\text{Error}}{1}{r} \lesssim
  \sum \mixednormlocal{\psi^{\pm}_{l} \left( \tfrac{M}{c^{2}} - 1 \right) 
  \overline{\psi^{\pm}_{l}}}{1}{r} + \tfrac{\sigma(T)}{c^{1/2}} 
  Y_{T}^{2}.
$$
The last term is certainly $o(1)$, and for the first term we write
$$
  \Lxpnorm{\psi^{\pm}_{l} \left( \tfrac{M}{c^{2}} - 1 \right) 
  \overline{\psi^{\pm}_{l}}}{r}
  \lesssim \Sobnorm{\psi^{\pm}}{1} \Sobnorm{\left( \tfrac{M}{c^{2}} - 1 \right) 
  \psi^{\pm}_{l}}{1}.
$$
Since the last factor is $\lesssim \Sobnorm{\psi^{\pm} - v^{\pm}}{1}
+ \Sobnorm{\left( \tfrac{M}{c^{2}} - 1 \right) v^{\pm}_{l}}{1}$, where we 
used Lemma \ref{LowHighLemma}\eqref{HL1}, it only remains to check
$$
  \energylocal{\left( \tfrac{M}{c^{2}} - 1 \right) v^{\pm}_{l}} = o(1),
$$
but this follows from the dominated convergence theorem.
This concludes the proof of \eqref{ReducedConvergence}, hence \eqref{Convergence}.

Then \eqref{DeltaA0minusu}, hence 
\eqref{A02uConvergence}, follows from a straightforward modification of the proof
of \eqref{Difference}, taking into account the estimates \eqref{errorBounds}
in the proof of Lemma \ref{A0Prop}, which hold for $1 \le r \le 3/2$. 
Thus, we can take $L_{t}^{\infty}$ of \eqref{Error}, instead of 
$L_{t}^{1}$.

This concludes the proof of Theorem \ref{Thm2}.

\section*{Appendix}

As mentioned in section \ref{MainResult}, the global existence
result of Klainerman and Machedon \cite{KM} was for the massless
KGM system.
In this appendix we show how their argument can be modified to handle the massive case.

First, the arguments relying on the conservation of energy require no change.
Thus, \cite[Proposition 1.1]{KM} holds as stated, and in fact the proof is easier
in the massive case, since now the KGM energy includes the $L^2$ norm of $\phi$.

The problem therefore reduces to proving local well-posedness for data with
$\mathcal I_0 < \infty$, where
$$
  \mathcal I_0 = \Sobdotnorm{\A\smallinit}{1} + \twonorm{\partial_t 
  \A\smallinit}{}
  + \Sobnorm{\phi\smallinit}{1} + \twonorm{\partial_t \phi\smallinit}{}.
$$
This is essentially what is proved in \cite[section 4]{KM}, and the argument
there is easily modified to handle the massive case. Let us give the details.
As in \cite{KM}, we set $c = 1$. Let $m > 0$ be the rest mass. Then
we have to add the linear term $m^2 \phi$ to the right hand side of \cite[Eq. (4.1b)]{KM},
which then corresponds to our equation \eqref{KGMc} (but with $c = 1$).
Then Propositions 4.1--4.4 in \cite{KM} hold as stated. The proofs only require
a few extra lines to treat the term $m^2 \phi$.

Consider Proposition 4.1. It is reduced to an inequality (see \cite[Eq. 
(4.3)]{KM}) which reads, in our notation,
$$
  f(T) \le \sigma(T) P (\mathcal I_0 + f(T) ),
$$
where $f(T) = \mixednormlocal{\square \A}{1}{2} + \mixednormlocal{\square \phi}{1}{2}$.
To extend this to the massive case, we only have to verify
$\mixednormlocal{m^2\phi}{1}{2} \le \sigma(T) P (\mathcal I_0 + f(T) )$.
To this end, write
$\mixednormlocal{m^2\phi}{1}{2} \le m^2 T \energylocal{\phi}$
and use the energy inequality
$$
  \energylocal{\phi}
  \lesssim \Sobnorm{\phi\smallinit}{1} + (1 + T) \twonorm{\partial_t \phi\smallinit}{}
  + (1 + T) \mixednormlocal{\square \phi}{1}{2}.
$$

Proposition 4.2 is a corresponding estimate for a difference of two solutions,
and since we have only added a linear term, the same changes apply there.
Finally, Propositions 4.3 and 4.4 are corollaries of Proposition 4.2.

\end{document}